\documentclass[12pt, a4paper]{amsart}

\usepackage[colorlinks, linkcolor=blue, anchorcolor=blue, citecolor=green]{hyperref}
\usepackage[normalem]{ulem}
\usepackage{graphics,epic}
\usepackage{amsmath,amssymb, amsthm,mathrsfs}
\usepackage[all,2cell]{xy}
\usepackage{shorttoc}
\usepackage{tikz}
\usetikzlibrary{matrix,positioning,decorations.markings,arrows,decorations.pathmorphing,	
	backgrounds,fit,positioning,shapes.symbols,chains,shadings,fadings,calc}
\tikzset{->-/.style={decoration={  markings,  mark=at position #1 with
			{\arrow{>}}},postaction={decorate}}}
\tikzset{-<-/.style={decoration={  markings,  mark=at position #1 with
			{\arrow{<}}},postaction={decorate}}}

\newtheorem{theorem}{Theorem}[section]
\newtheorem*{theorem*}{Theorem}

\newtheorem*{conjecture*}{Conjecture}

\newtheorem{thm}[theorem]{Theorem}
\newtheorem{lem}[theorem]{Lemma}

\newtheorem{df}[theorem]{Definition}

\newcommand{\ie}{{\em i.e.}\ }

\newcommand{\opname}[1]{\operatorname{\mathsf{#1}}}

\newcommand{\Arc}{\opname{Arcs}}
\newcommand{\Arcseg}{\opname{Arcs}}

\newcommand{\Z}{\mathbb{Z}}

\renewcommand{\P}{\mathbb{P}}

\newcommand{\fT}{\mathbf{T}}
\newcommand{\fS}{\mathbf{S}}
\newcommand{\fM}{\mathbf{M}}
\newcommand{\TT}{\mathbf{T}}
\renewcommand{\SS}{\mathbf{S}}
\newcommand{\MM}{\mathbf{M}}
\newcommand{\ba}{\mathbf{a}}
\newcommand{\bl}{\mathbf{l}}

\newcommand{\Int}{\opname{Int}}
\newcommand{\dInt}{\opname{Int}^{\circ}}
\newcommand{\dIntt}{\opname{Int}}
\newcommand{\Intv}{\underline{\opname{Int}}}
\newcommand{\dIntv}{\underline{\opname{Int}}^{\circ}}
\newcommand{\dInttv}{\underline{\opname{Int}}}

\newcommand{\den}{\opname{den}}

\newcommand{\cm}{{\mathcal M}}
\newcommand{\cn}{{\mathcal N}}
\newcommand{\cM}{{\mathcal M}}
\newcommand{\CM}{{\mathcal M}}
\newcommand{\cN}{{\mathcal N}}

\newcommand{\ma}{\mathcal{A}}

\newcommand{\nn}{node[black]{$\bullet$}}

\makeatletter
\newcommand*\bigcdot{\mathpalette\bigcdot@{.5}}
\newcommand*\bigcdot@[2]{\mathbin{\vcenter{\hbox{\scalebox{#2}{$\m@th#1\bullet$}}}}}
\makeatother

\begin{document}

\title{Denominator conjecture for some surface cluster algebras}

\author{Changjian Fu}
\address{Changjian Fu\\Department of Mathematics\\SiChuan University\\610064 Chengdu\\P.R.China}
\email{changjianfu@scu.edu.cn}
\author{Shengfei Geng}
\address{Shengfei Geng\\Department of Mathematics\\SiChuan University\\610064 Chengdu\\P.R.China}
\email{genshengfei@scu.edu.cn}
\maketitle

\begin{abstract}
The denominator conjecture, proposed by Fomin and Zelevinsky, says that for a cluster algebra, the cluster monomials are uniquely determined by their denominator vectors with respect to an initial cluster. In this paper, for a cluster algebra from a marked surface with at least three boundary marked points, we establish this conjecture with respect to a given strong admissible tagged triangulation.

 
\end{abstract}
\tableofcontents
\section{Introduction} 
Cluster algebra was introduced by Fomin and Zelevinsky \cite{FZ02}   in
 order to give an algebraic
 framework of canonical basis and positivity. 
It is
in essence a commutative ring with a distinguished countable
family of generators, called {\em cluster variables}, where each cluster variable can be obtained from the initial cluster variables by a sequences of mutations. The cluster
variables can be grouped into overlapping sets of equal finite
size $n$, called {\em clusters}.
Let $\mathcal{A}$ be a cluster algebra with trivial coefficients and  $\mathbf{x}=(x_1,\dots,x_n)$ a cluster of $\mathcal{A}$. 
 It follows from the Laurent phenomenon  \cite{FZ02} that every cluster variable $y\in \ma$   can be uniquely written as $$y=\frac{f(x_1,\dots,x_n)}{\prod_{i=1}^nx_i^{d_i}}, $$ where $f(x_1,\dots,x_n)\in \mathbb{Z}[x_1,\dots,x_n]$ can not be divisible by any $x_i$ and $d_i\in \mathbb Z$ for any $1\leq i\leq n$.
The vector $$\den(y)=(d_1,\cdots,d_n)$$ is called the {\em denominator vector} of $y$ with respect to $\mathbf{x}$. A {\it cluster monomial} of $\mathcal{A}$ is a monomial of cluster variables belonging to the same cluster.
Let $(y_1,\dots,y_n)$ be a cluster and $\mathbf{m}={\prod_{i=1}^ny_i^{a_i}}$  a cluster monomial, where $a_1,\dots, a_n$ are  non-negative integer.  The {\em denominator vector} of $\mathbf{m}$ with respect to $\mathbf{x}$ is defined as
$$\den(\mathbf{m})=\Sigma_{i=1}^na_i\den(y_i).$$
Fomin and Zelevinsky proposed the following denominator conjecture:
\begin{conjecture*}\cite[Conjecture 4.17]{FZ04} \cite[Conjecture 7.6]{FZ07} Let $\mathcal{A}$ be a cluster algebra and $\mathbf{x}$  an arbitrary cluster. Then different cluster monomials have different denominator vectors with respect to $\mathbf{x}$.
\end{conjecture*}
The denominator conjecture remains a challenging and open problem in the field of cluster algebras.
 For acyclic cluster algebras with acyclic initial clusters, the conjecture has been confirmed. This includes verification for cluster algebras of rank 2 by Sherman and Zelevinsky \cite{SZ}, for acyclic skew-symmetric cluster algebras with arbitrary acyclic clusters by Caldero and Keller \cite{CK06,CK08}, and for cluster algebras of finite type with bipartite initial clusters by Fomin and Zelevinsky \cite{FZ07}.
By using the quantum cluster character, 
Rupel and Stella \cite{RS} further removed the condition of skew-symmetric type and verified the conjecture for all acyclic cluster algebras with respect to acyclic clusters. 
In a recent development, the authors  \cite{FG2022} verified the denominator conjecture  for cluster algebras of types $\mathbb{A}$, $\mathbb{B}$, and $\mathbb{C}$ with respect to any initial cluster.  Additionally, a weaker version of the conjecture has been explored in several works, see \cite{GP, NS, FG, GY20, FGL21a} for instance.

Cluster algebras from marked surfaces are a special class of cluster algebras that were developed in \cite{FoG06,FoG09,FST08,FT18,GSV}. These algebras provide a bridge between the fields of algebra, combinatorics, and geometry, and they have been a subject of intense study due to their rich structure and the surprising connections they revealed, see  \cite{FST08, FT18, FeST12,  IL12, L09, MSW11, MSW13, M17} for instance.
A marked surface $(\SS,\MM)$ is a compact connected Reimann surface with a finite set of  boundary marked points and some interior marked points  called {\em punctures}. A {\it tagged triangulation} $\TT$ of $(\SS,\MM)$ is a maximal set of pairwise compatible tagged arcs. Each tagged triangulation $\TT$ of $(\SS,\MM)$ determines a cluster algebra $\ma_\TT$. If $\SS$ has non-empty boundary, there is a bijection between the set of cluster monomials  of $\ma_\TT$ and  the set of finite multisets of pairwise compatible tagged arcs over $(\SS,\MM,\TT)$.
 Moreover, the denominator vector of a cluster variable  coincides with the intersection vector of the corresponding tagged arc. More generally, the denominator vector of a cluster monomial  coincides with the intersection vector of the corresponding multiset of pairwise compatible tagged arcs, see \cite{FST08}.

 The aim of this paper is to establish the denominator conjecture for some cluster algebras arsing from marked surface with respect to a particular choice of initial seed.
 Let $(\SS,\MM)$ be a marked surface with boundary such that each boundary component  has  at least one marked point. Let $\TT$ be an admissible ideal triangulation of  $(\SS,\MM)$. 
For a finite multiset $\cm$ of pairwise compatible arcs such hat $\cm\cap \TT=\emptyset$,
 we   denote by $\Arcseg^1(\cM)$ the multiset of  all the  irreducible arc segments of arcs in $\cM$. When $\SS$ has no punctures, it has been proved  that $\cm$ is uniquely determined by  $\Arcseg^1(\cM)$ \cite{FG2022}. 
 Our first main result is a generalization of  the aforementioned result.
\begin{thm}\label{t:local to entire in background}(Theorem~\ref{t:local to entire}) Let $(\fS,\MM)$ be a marked surface equipped with an admissible ideal triangulation $\TT$. 
Let $\cm$ and $\cn$ be two finite multisets of pairwise compatible arcs such that  $\cm\cap\TT=\emptyset$ and $\cn\cap \TT=\emptyset$.  Then  $\cm=\cn$ if and only if $\Arcseg^1(\cM)=\Arcseg^1(\cN)$.
\end{thm}

 A tagged triangulation  $\TT$ of $(\SS,\MM)$  is called  {\it  strong admissible} if $\TT$ contains no loops and each puncture of $(\SS,\MM)$ is connected to a boundary marked point by a pair of conjugate arcs in $\TT$.
Then, we have 
\begin{thm}(Theorem~\ref{t:main theorem-admissible case}+ Theorem~\ref{t:main theorem with strong admissible initial seed})
Let $(\SS,\MM)$ be a marked surface with a strong admissible tagged triangulation $\TT$. 
\begin{itemize}
    \item [(1)] Let $\cm$ and $\cn$ be two finite multisets of pairwise compatible tagged arcs. If $\cm$ and $\cn$ have the same intersection vectors, then $\cm=\cn$.
    \item[(2)] Let $\ma_\TT$ be the cluster algebra corresponding to $(\SS,\MM,\TT)$. Then different cluster monomials of $\ma_\TT$ have different denominator vectors with respect to the initial seed defined by $\TT$.
\end{itemize}
\end{thm}
If $\SS$ is of genus $0$ with at least two boundary marked points, 
 one can deduce that such a strong admissible tagged triangulation  always exists. In the general case,
if $\SS$ possesses at least three boundary marked points, the existence of a strong admissible tagged triangulation is guaranteed (see Lemma~\ref{l:exsitence of strong admissible}).  Therefore, for a cluster algebra arising from a marked surface with at least three boundary marked points or from a genus $0$ marked surface with at least two boundary marked points, the denominator conjecture holds true for any given strong admissible tagged triangulation.

\vspace{0.2cm}
The structure of the paper is as follows. In Section~\ref{s:marked surfaces}, we review some fundamental definitions of marked surfaces. Specifically, we define a strong admissible ideal triangulation of a marked surface and discuss its existence. Additionally, we prove Theorem~\ref{t:local to entire in background}. In Section~\ref{s:Multiset of untagged arcs determined by their dimension vectors}, we demonstrate that for a marked surface with a strong admissible ideal triangulation, each multiset of pairwise compatible arcs is uniquely determined by its intersection vector. In Section~\ref{s:Multiset of tagged arcs determined by their intersection vectors}, we first revisit some basic definitions of tagged arcs, (strong admissible) tagged triangulations, and intersection vectors of tagged arcs. We then explore their connection with the untagged version. Finally, we show that for a marked surface with a strong admissible tagged triangulation, each multiset of pairwise compatible tagged arcs is uniquely determined by its intersection vector (Theorem \ref{t:main theorem-admissible case}). In Section~\ref{s:Denominator conjecture  for some cluster algebras}, we establish the denominator conjecture for  cluster algebras arising from  marked surfaces with  strong admissible tagged triangulations (Thoerem \ref{t:main theorem with strong admissible initial seed}).


\section{Marked surfaces}\label{s:marked surfaces}
\subsection{Marked surface}
Let $\SS$ be a connected oriented (2-dimensional) Riemann surface with boundary. Fix a finite set $\MM$ of marked points on $\SS$ that includes at least one marked point on each boundary component, plus possibly some interior points, called {\em punctures}.
As \cite{FST08}, we assume $(\bold{S}, \bold{M})$ is none of the following:
\begin{itemize}
\item an unpunctured or once-punctured monogon;
\item an unpunctured digon;
\item an unpunctured triangle.
(In this paper, a $m$-gon is a disk with $m$ marked points on the boundary.)
\end{itemize}

\begin{df}
A (simple) arc $\gamma$ in $(\SS, \MM)$ is a curve $\gamma$ in $\bold{S}$ such that
\begin{itemize}
 \item the endpoints of $\gamma$ are marked points in $\MM$;
\item $\gamma$ does not intersect itself, except that its endpoints may coincide;
\item except for the endpoints, $\gamma$ is disjoint from $\bold{M}$ and from the boundary of $\bold{S}$;
\item  $\gamma$ does not cut out an unpunctured monogon or an unpunctured digon. (In other
words, $\gamma$ is not contractible into $\bold{M}$ or onto the boundary of $\SS$.)
\end{itemize}
\end{df}
An arc $\gamma$ is a {\it loop} if its endpoints coincide, where  the common endpoint is known as the  {\it base point} of $\gamma$.

\subsection{Admissible ideal triangulation}
Two arcs in $(\bold{S}, \bold{M})$ are called {\em compatible} if they do not intersect in the interior of $\bold{S}$, more precisely, there are curves in their respective isotopy classes which do not intersect in the interior of $\bold{S}$. It is clear that  each arc is compatible with itself. An  {\it ideal triangulation $\TT$} of $(\SS,\MM)$ is a maximal set of pairwise compatible arcs.

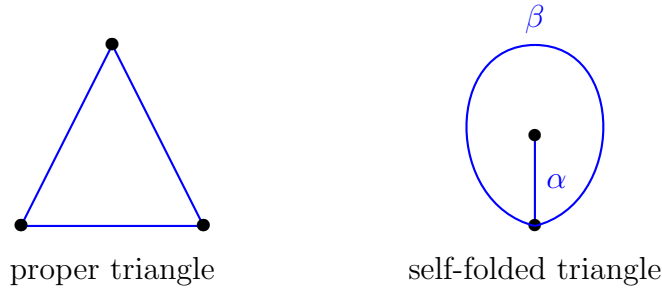
\begin{figure}[h]
\begin{minipage}[t]{0.3\linewidth} 

\begin{tikzpicture}[xscale=0.6,yscale=0.6]

			\draw[blue,thick] (-2,-2)\nn to (2,-2)\nn;
			
		\draw[blue,thick] (-2,-2)\nn to (0,2)\nn;
			\draw[blue,thick] (2,-2)\nn to (0,2)\nn;

\node at (0,-3){proper triangle};
	\end{tikzpicture}
\end{minipage}%
\begin{minipage}[t]{0.3\linewidth} 
\centering
\begin{tikzpicture}[xscale=0.6,yscale=0.6]
			
			\draw[blue,thick] (0,-2)\nn to (0,0)\nn;
			
\draw[thick,color=blue] (0,-2) .. controls (-2,-1.5) and (-2,2) .. (0,2);	
\draw[thick,color=blue] (0,-2) .. controls (2,-1.5) and (2,2) .. (0,2);		
\draw[thick,blue](0,-1)node[right]{$\alpha$}(0,2)node[above]{$\beta$};

\node at (0,-3){self-folded triangle};

	\end{tikzpicture}
\end{minipage}

\caption{Basic triangles}\label{f:triangle}
\end{figure}

A triangle in $\TT$ has three distinct sides unless it is a self-folded triangle as Figure~\ref{f:triangle}, where we call $\alpha$ the folded side  and $\beta$ the remaining side.
An ideal triangulation of $(\bold{S}, \bold{M})$ is called {\it admissible} if each puncture of $(\bold{S}, \bold{M})$ is enclosed in a self-folded triangle.

\subsection{Strong admissible ideal triangulation}
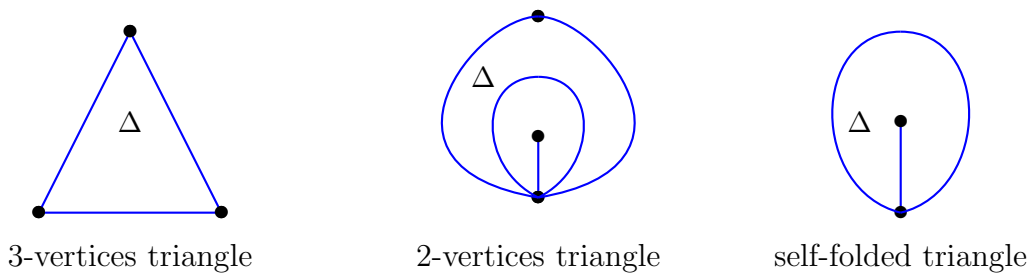
\begin{figure}[h]
\begin{minipage}[t]{0.3\linewidth} 

\begin{tikzpicture}[xscale=0.6,yscale=0.6]

			\draw[blue,thick] (-2,-2)\nn to (2,-2)\nn;
			
		\draw[blue,thick] (-2,-2)\nn to (0,2)\nn;
			\draw[blue,thick] (2,-2)\nn to (0,2)\nn;

\draw[black](0,0)node{{\small$\Delta$}};	

\node at (0,-3){3-vertices triangle};
	\end{tikzpicture}

\end{minipage}%
\begin{minipage}[t]{0.3\linewidth} 

\begin{tikzpicture}[xscale=0.8,yscale=0.8]
		
			\draw[blue,thick] (0,-2)\nn to (0,-2)\nn;
				\draw[blue,thick] (0,1)\nn to (0,1)\nn;
	\draw[blue,thick] (0,-1)\nn to (0,-2)\nn;

\draw[thick,color=blue] (0,-2) .. controls (-1,-1.5) and (-1,0) .. (0,0);	
\draw[thick,color=blue] (0,-2) .. controls (1,-1.5) and (1,0) .. (0,0);

\draw[thick,color=blue] (0,-2) .. controls (-3,-1.5) and (-1,1) .. (0,1);	
\draw[thick,color=blue] (0,-2) .. controls (3,-1.5) and (1,1) .. (0,1);	

\draw[black](-0.9,0)node{{\small$\Delta$}};	

\node at (0,-3){2-vertices triangle};

		\end{tikzpicture}
		
\end{minipage}
\begin{minipage}[t]{0.3\linewidth} 
\centering
\begin{tikzpicture}[xscale=0.6,yscale=0.6]
			
			\draw[blue,thick] (0,-2)\nn to (0,0)\nn;
			
\draw[thick,color=blue] (0,-2) .. controls (-2,-1.5) and (-2,2) .. (0,2);	
\draw[thick,color=blue] (0,-2) .. controls (2,-1.5) and (2,2) .. (0,2);		

\draw[black](-0.9,0)node{{\small$\Delta$}};	
\node at (0,-3){self-folded triangle};

	\end{tikzpicture}
\end{minipage}

\caption{Basic triangle under strong admissible triangulation}\label{f:Basic triangle under strong admissible triangulation}
\end{figure}

A loop that serves as the remaining side of a self-folded triangle is labeled as a {\it type $\mathbf{I}$} loop. Conversely, a loop that does not function as the remaining side of any self-folded triangle is termed a {\it type $\mathbf{II}$} loop.
 An  admissible ideal triangulation $\TT$ is called {\it strong} if each loop in $\TT$ is  a type $\mathbf{I}$  loop. So if $\TT$ is a strong admissible triangulation,  it is easy to get that each basic triangle is either a proper triangle with three different vertices (simply called 3-vertices triangle), a self-folded triangle or a proper triangle close wrapping  some self-folded triangle (simply called 2-vertices triangle)(cf. Figure~\ref{f:Basic triangle under strong admissible triangulation}).

By the assumption on $(\SS,\MM)$,  it is clear that the admissible ideal triangulation always exists. However, the strong admissible ideal triangulation may not exist.  But if $\SS$ is genus 0 with at least two boundary points, it is straightforward to obtain a strong admissible  ideal triangulation. In the general case, we have
\begin{lem} \label{l:exsitence of strong admissible}
Let $(\SS,\MM)$ be a connected marked surface with at least three marked points on the boundary, then there is a strong admissible ideal triangulation $\TT$ on $(\SS,\MM)$.
\end{lem}
\begin{proof}
 Let $\TT_1$   be an admissible triangulation on $(\SS,\MM)$, \ie a  triangulation where each puncture is enclosed within a self-folded triangle.
 It is straightforward to obtain such a triangulation.
 Let $l$ represent the count  of type $\mathbf{II}$   loops in $\TT_1$.  If $l=0$, $\TT_1$ is just we wanted.
Now assume $l>0$. By induction, it is  sufficient to demonstrate that there exists an admissible ideal triangulation $\TT_2$  with fewer type $\mathbf{II}$  loops than $l$.
   
   Take $\alpha$ as a type $\mathbf{II}$   loop  in $\TT_1$. The arc $\alpha$ forms an edge of two  proper triangles, $\Delta_1$ and  $\Delta_2$. Let the base point of $\alpha$
   be $A$. Label  the other vertex of $\Delta_1$ as $B$, the other vertex of $\Delta_2$ as $C$. We then encounter only four possible cases:

\begin{itemize}
    \item Case 1: $B\neq A$, $C\neq A$ and $B\neq C$ (cf. Figure \ref{f:Cases of loop}: Case 1).
     \item Case 2: $B=A\neq C$ or $C=A\neq B$ (cf. Figure \ref{f:Cases of loop}: Case 2).
 
     \item Case 3: $B=C=A$ (cf. Figure \ref{f:Cases of loop}: Case 3). 
      \item Case 4: $B=C\neq A$ (cf. Figure \ref{f:Cases of loop}: Case 4).
      
\end{itemize}
     
\begin{figure}[htpb]
   \begin{minipage}{0.23\linewidth}
		
\begin{tikzpicture}[xscale=0.6,yscale=0.6]
	
			\draw[thick] (0,0)\nn to (0,4)\nn;
   \draw[thick] (-2,2)\nn to (0,0)\nn;
   \draw[thick] (2,2)\nn to (0,0)\nn;
      \draw[thick] (-2,2)\nn to (0,4)\nn;
   \draw[thick] (2,2)\nn to (0,4)\nn;
\draw[blue](0,0)node[below]{$A$} (2,2)node[right]{$C$}(-2,2)node[left]{$B$}(0,4)node[above]{$A$}(0,2)node[right]{$\alpha$};

\draw[blue](0,-1)node[below]{Case 1}; 
		\end{tikzpicture}

\end{minipage}		
 \begin{minipage}{0.23\linewidth}
		
\begin{tikzpicture}[xscale=0.6,yscale=0.6]
	
			\draw[thick] (0,0)\nn to (0,4)\nn;
   \draw[thick] (-2,2)\nn to (0,0)\nn;
   \draw[thick] (2,2)\nn to (0,0)\nn;
      \draw[thick] (-2,2)\nn to (0,4)\nn;
   \draw[thick] (2,2)\nn to (0,4)\nn;
\draw[blue](0,0)node[below]{$A$} (2,2)node[right]{$C$}(-2,2)node[left]{$A$}(0,4)node[above]{$A$}(0,2)node[right]{$\alpha$};

\draw[blue](0,-1)node[below]{Case 2}; 
		\end{tikzpicture}
		

\end{minipage}	
 \begin{minipage}{0.23\linewidth}
		
\begin{tikzpicture}[xscale=0.6,yscale=0.6]
	
			\draw[thick] (0,0)\nn to (0,4)\nn;
   \draw[thick] (-2,2)\nn to (0,0)\nn;
   \draw[thick] (2,2)\nn to (0,0)\nn;
      \draw[thick] (-2,2)\nn to (0,4)\nn;
   \draw[thick] (2,2)\nn to (0,4)\nn;
\draw[blue](0,0)node[below]{$A$} (2,2)node[right]{$A$}(-2,2)node[left]{$A$}(0,4)node[above]{$A$}(0,2)node[right]{$\alpha$};

\draw[blue](0,-1)node[below]{Case 3}; 
		\end{tikzpicture}
		

\end{minipage}	
 \begin{minipage}{0.23\linewidth}
		
\begin{tikzpicture}[xscale=0.6,yscale=0.6]
	
			\draw[thick] (0,0)\nn to (0,4)\nn;
   \draw[thick] (-2,2)\nn to (0,0)\nn;
   \draw[thick] (2,2)\nn to (0,0)\nn;
      \draw[thick] (-2,2)\nn to (0,4)\nn;
   \draw[thick] (2,2)\nn to (0,4)\nn;
\draw[blue](0,0)node[below]{$A$} (2,2)node[right]{$B$}(-2,2)node[left]{$B$}(0,4)node[above]{$A$}(0,2)node[right]{$\alpha$};

\draw[blue](0,-1)node[below]{Case 4}; 
		\end{tikzpicture}
		

\end{minipage}	

	\caption{Cases of $\alpha$}\label{f:Cases of loop}
\end{figure}
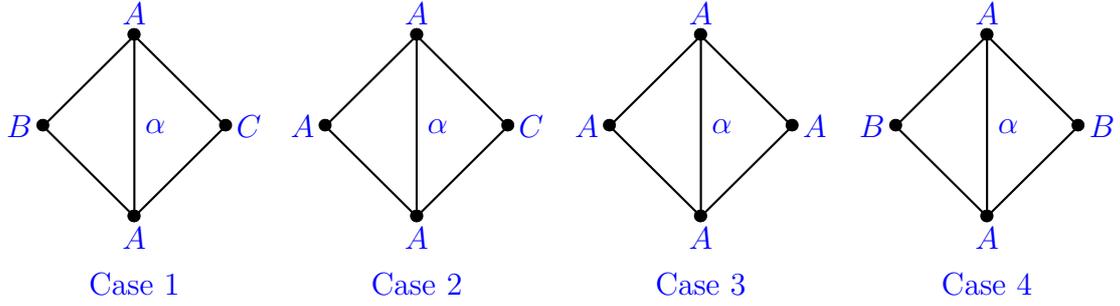

Every type $\mathbf{II}$  loop behaves like $\alpha$, fitting exactly into one of the four cases mentioned. If there is a type $\mathbf{II}$  loop similar to $\alpha$ in Case 1 or Case 2, we can, without any loss of generality, assume this loop to be $\alpha$ itself. In this instance, we simply substitute $\alpha$ with the alternative diagonal $\alpha'$. This action results in a new admissible triangulation $\TT_2$. Since $\alpha'$ is not a loop, $\TT_2$ will have fewer type $\mathbf{II}$  loops than $\TT_1$.

Now, let us assume that there are no type $\mathbf{II}$  loops behaving like $\alpha$ in Case 1 or Case 2. We then claim that  a type $\mathbf{II}$  loop like $\alpha$ cannot exist in Case 3  either. This is because $\SS$ features at least three marked points on its boundary, meaning there is at least one arc in $\TT_1$ that is neither a loop nor a radius of a self-folded triangle. Given that $\SS$ is dissected into a finite number of triangles by the arcs in $\TT_1$ and $\SS$ is connected, it follows that there must be a triangle in Case 3, with at least one edge fitting into Case 2, which is a contradiction.

Considering this, we can proceed with the understanding that all type $\mathbf{II}$  loops in $\TT_1$ are in Case 4. Because $\SS$ has at least three marked points on its boundary and $\SS$ is connected, there is a triangle $\Delta$ with three distinct vertices and there is an edge $\beta$ of $\Delta$ that is part of another triangle $\Delta'$ such that $\Delta'$ includes an edge $\alpha$ which is a type $\mathbf{II}$  loop. By assumption, $\alpha$ is in Case 4, if we flip $\beta$ and then flip $\alpha$, we eliminate one type $\mathbf{II}$ loop without creating any additional type $\mathbf{II}$  loops (cf. Figure~\ref{f:Flip in Case 4}). Consequently, we can define $\TT_2$ as the triangulation obtained from $\TT_1$ by first flipping $\beta$ and then flipping $\alpha$. Hence, $\TT_2$ is a new admissible triangulation of $(\SS,\TT)$ with 
fewer type $\mathbf{II}$  loops than $\TT_1$.
\begin{figure}[htpb]
   \begin{minipage}{0.3\linewidth}
		
\begin{tikzpicture}[xscale=0.6,yscale=0.6]
	
			\draw[thick] (0,0)\nn to (0,4)\nn;
   \draw[thick] (-2,2)\nn to (0,0)\nn;
   \draw[thick] (2,2)\nn to (0,0)\nn;
      \draw[thick] (-2,2)\nn to (0,4)\nn;
   \draw[thick] (2,2)\nn to (0,4)\nn;
\draw[blue](0,0)node[below]{$A$} (2,2)node[below]{$B$}(-2,2)node[below]{$B$}(0,4)node[above]{$A$}(0,2)node[right]{$\alpha$};
\draw[thick] (-2,4)\nn to (0,4)\nn;
 \draw[thick] (-2,4)\nn to (-2,2)\nn;
 \draw[blue](-2,4)node[above]{$C$}(-1,3)node[below]{$\beta$};
  \draw[red](3,2)node[right]{$\longrightarrow$}(4,2)node[above]{flip $\beta$};
 \end{tikzpicture}

\end{minipage}		
 \begin{minipage}{0.3\linewidth}
		
\begin{tikzpicture}[xscale=0.6,yscale=0.6]
	
	\draw[thick] (0,0)\nn to (0,4)\nn;
   \draw[thick] (-2,2)\nn to (0,0)\nn;
   \draw[thick] (2,2)\nn to (0,0)\nn;
      \draw[thick] (-2,4)\nn to (0,0)\nn;
   \draw[thick] (2,2)\nn to (0,4)\nn;
\draw[blue](0,0)node[below]{$A$} (2,2)node[below]{$B$}(-2,2)node[below]{$B$}(0,4)node[above]{$A$}(0,2)node[right]{$\alpha$};
\draw[thick] (-2,4)\nn to (0,4)\nn;
 \draw[thick] (-2,4)\nn to (-2,2)\nn;
 \draw[blue](-2,4)node[above]{$C$};
  \draw[red](3,2)node[right]{$\longrightarrow$}(4,2)node[above]{flip $\alpha$};
		\end{tikzpicture}
		

\end{minipage}	
 \begin{minipage}{0.3\linewidth}
		
\begin{tikzpicture}[xscale=0.6,yscale=0.6]
	
				\draw[thick] (2,2)\nn to (-2,4)\nn;
   \draw[thick] (-2,2)\nn to (0,0)\nn;
   \draw[thick] (2,2)\nn to (0,0)\nn;
      \draw[thick] (-2,4)\nn to (0,0)\nn;
   \draw[thick] (2,2)\nn to (0,4)\nn;
\draw[blue](0,0)node[below]{$A$} (2,2)node[below]{$B$}(-2,2)node[below]{$B$}(0,4)node[above]{$A$}(0,2)node[right]{$\alpha$};
\draw[thick] (-2,4)\nn to (0,4)\nn;
 \draw[thick] (-2,4)\nn to (-2,2)\nn;
 \draw[blue](-2,4)node[above]{$C$};
		\end{tikzpicture}

\end{minipage}	

	\caption{Flip in Case 4}\label{f:Flip in Case 4}
\end{figure}
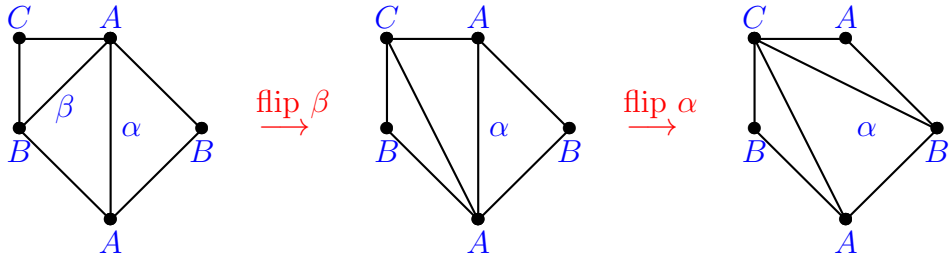
\end{proof}

\subsection{Arc segments}
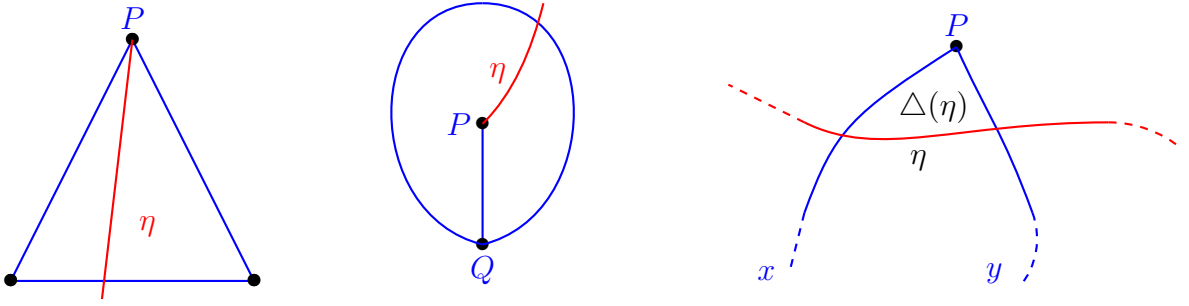
\begin{figure}[h]

	\begin{minipage}{0.3\linewidth} 
		\begin{tikzpicture}[xscale=0.8,yscale=0.8]
			
			\draw[blue,thick] (-2,-2)\nn to (2,-2)\nn;
			
		\draw[blue,thick] (-2,-2)\nn to (0,2)\nn;
			\draw[blue,thick] (2,-2)\nn to (0,2)\nn;

		\draw[blue](0,2)node[above]{$P$};

  \draw[red](0.25,-0.75) node[below]{$\eta$};

	\draw[red,thick] (0,2) to (-0.5,-2.3);

		\end{tikzpicture}
		\end{minipage} 
\begin{minipage}{0.3\linewidth}
\begin{tikzpicture}[xscale=0.8,yscale=0.8]
			\draw[blue,thick] (0,-2)\nn to (0,0)\nn;
			
\draw[thick,color=blue] (0,-2) .. controls (-2,-1.5) and (-2,2) .. (0,2);	
\draw[thick,color=blue] (0,-2) .. controls (2,-1.5) and (2,2) .. (0,2);		
\draw[thick,blue](0,-2)node[below]{$Q$}(0,0)node[left]{$P$};
\draw[red](0.25,0.5)node[above]{$\eta$};

\draw[thick,color=red] (0,0) .. controls (0.5,0.5) and (0.8,1.2) .. (1,2);	

	\end{tikzpicture}
\end{minipage}
\begin{minipage}{0.3\linewidth}
    \begin{tikzpicture}
			\draw[thick, color=blue](10,0)\nn node[above]{$P$};	
    
    \draw[thick,color=blue] (10,0) .. controls (8.6,-0.9) and (8.4,-1.1) .. (8,-2.2);
				\draw[thick,color=blue,dashed] (8,-2.2) -- (7.8,-3);
				\draw[thick,color=blue] (10,0) .. controls (10.4,-0.9) and (10.6,-1.1) .. (11,-2.2);
				\draw[thick,color=blue,dashed] (11,-2.2) arc (20:-40:1);
				\draw[thick,color=red] (8,-1) .. controls (9,-1.5) and (10,-1) .. (12,-1);
				\draw[thick,color=red,dashed] (12,-1) arc (90:50:1.5);
				\draw[thick,color=red,dashed] (8,-1) -- (7,-0.5);
				\node at (9.5, -1.5){$\eta$};
				\node at (7.5,-3){{\color{blue}$x$}};
				\node at (10.5, -3){{\color{blue}$y$}};
				\draw (9.7,-.8)node{$\triangle(\eta)$};
			\end{tikzpicture}
\end{minipage}
  \caption{Irreducible arc segments}\label{fig:irreducible arc segments}
  \end{figure}

Let $\fT$ be an ideal triangulation of $(\fS,\fM)$. Every arc is decomposed by $\fT$ into segments. The equivalence of arcs induces an equivalence relation on  segments lying in the same region of $\fT$. An equivalence class of a segment is called an {\em irreducible arc segment} (with respect to $\fT$).
In Figure \ref{fig:irreducible arc segments}, we have list all the possible irreducible segments, where the first two segments appear at the ends of arcs, while the third one appear in the middle of an arc.

Note that if 
the endpoints of an irreducible arc segment $\eta$ are in the interiors of non-boundary edges $x,y$  of a triangle $\Delta$, then $x,y$ have a common endpoint $p_\eta:=P\in \MM$ and $\eta$ cuts out an angle from $\Delta$ (cf. Figure \ref{fig:irreducible arc segments}). We denote by $\triangle(\eta)$ the local triangle cut out by $\eta$ and denoted by $\angle \eta$ the angle in $\Delta(\eta)$ with vertex $p_\eta$ and two edges $x, y$. 	

For a finite multiset $\cm$ of arcs that are pairwise compatible,
 we   denote by $\Arcseg^1(\cM)$ the multiset of  all the  irreducible arc segments of arcs in $\cM$.

\begin{thm}\label{t:local to entire} Consider a marked surface $(\fS,\MM)$ equipped with an admissible ideal triangulation $\TT$. 
Let $\cm$, $\cn$ be two finite multisets,  where each one is consisted of pairwise compatible arcs without element from $\TT$.  Then  $\cm=\cn$ if and only if $\Arcseg^1(\cM)=\Arcseg^1(\cN)$.
\end{thm}
 \begin{proof}
 Given that $\cm=\cn$, it is evident that $\Arcseg^1(\cM)=\Arcseg^1(\cN)$.
 Now, assume the converse, that is,  $\Arcseg^1(\cM)=\Arcseg^1(\cN)$.
Let $P_1,P_2,\cdots,P_t$ be all the punctures of $\SS$. If $t=0$, then the result follows from  \cite[Lemma2.8]{FG2022}.

   Assume that $t>0$. Since  $\TT$ is admissible, each puncture is contained in a self-folded triangle.  For each $1\leq k\leq t$,   let us identify the self-folded triangle enclosing $P_k$ as $\Delta_k$, with  $\gamma_k$ as the radius and $l_k$ as  the remaining side.
   Denote by $Q_k$ the base point of $l_k$.
   Let $\delta_1,\cdots,\delta_j$ be all the irreducible arc segments in  $\Arcseg^1(\cM)=\Arcseg^1(\cN)$ with one endpoint $P_k$.   Note that for each $1\leq i\leq j$, the connection  of each $\delta_i$ and $\gamma_k$, denoted by $\delta_i\circ\gamma_k$,
    divides $\Delta_k$ into two parts. Fix a direction of $\delta_i\circ\gamma_k$ with $Q_k$ as the starting point of $\gamma_k\circ \delta_i$, denoted by $L_i$ (resp. $R_i$) the left part (resp. right part)  of $\Delta_k$ (cf. Figure~\ref{f:Some changes in self-folded triangle}).
 We then replace  each puncture $P_k$ by a boundary component $B_k$ such that $B_k$ lies in the intersection  $\bigcap_{i=1}^jL_i$,  featuring a solitary marked point, still referred to as $P_k$, while maintaining  $\gamma_k$ and $l_k$ unchanged (cf. Figure~\ref{f:Some changes in self-folded triangle}.) 
   Consequently, $(\fS,\MM,\TT)$ transfer into a new marked surface $(\fS',\MM',\TT')$  with a tiling $\TT'$, where $\SS'$ has no punctures. 
 It is clear that $\delta_i$ is permissible (cf. ~\cite[2.2]{FG2022}) in $(\fS',\MM',\TT')$ for each $1\leq i\leq j$. Then in $(\fS',\MM',\TT')$,  the arcs in $\cM$ and $\cN$ are permissible (cf. ~\cite[Definition 2.2]{FG2022}). It is straightforward to observe that $\Arcseg^1(\cM)=\Arcseg^1(\cN)$ in $(\SS',\MM',\TT')$. Since $\SS'$ has no punctures,  by ~\cite[Lemma 2.8]{FG2022}, we conclude that $\cm=\cn$ in $(\SS',\MM',\TT')$, and thus $\cm=\cn$ in $(\SS,\MM,\TT)$.
\end{proof}

\begin{figure}[h]
\begin{minipage}[t]{0.3\linewidth} 

\begin{tikzpicture}[xscale=0.8,yscale=0.8]

			\draw[blue,thick,->-=.5,>=stealth] (0,-2)\nn to (0,0)\nn;
			
\draw[thick,color=blue] (0,-2) .. controls (-2,-1.5) and (-2,2) .. (0,2);	
\draw[thick,color=blue] (0,-2) .. controls (2,-1.5) and (2,2) .. (0,2);		
\draw[thick,blue](0,-1)node[left]{$\gamma_k$}(0,2)node[above]{$l_k$}(0,0)node[left]{$P_k$}(0,-2)node[below]{$Q_k$}(1,0)node[below]{$\delta_i$}(-0.3,1)node[left]{$L_i$}(0.2,-1)node[right]{$R_i$};

\draw[red,thick,->-=.5,>=stealth] (0,0)\nn to (2,0);

	\end{tikzpicture}
\end{minipage}%
\begin{minipage}[t]{0.2\linewidth} 
\begin{tikzpicture}[xscale=0.8,yscale=0.8]
\draw[white,thick,=>] (0,2) to (2,2);
\draw[red,thick,->] (0,0) to (2,0);
\draw[white,thick,->] (0,-2) to (2,-2);

	\end{tikzpicture}
\end{minipage}%
\begin{minipage}[t]{0.3\linewidth} 
\begin{tikzpicture}[xscale=0.8,yscale=0.8]
			
	\draw[thick,fill=black!20] (0,0.45) circle(12pt);

			\draw[blue,thick,->-=.5,>=stealth] (0,-2)\nn to (0,0)\nn;
   
\draw[thick,color=blue] (0,-2) .. controls (-2,-1.5) and (-2,2) .. (0,2);	
\draw[thick,color=blue] (0,-2) .. controls (2,-1.5) and (2,2) .. (0,2);		
\draw[thick,blue](0,-1)node[left]{$\gamma_k$}(0,2)node[above]{$l_k$}(0,0)node[left]{$P_k$}(0,-2)node[below]{$Q_k$}(1,0)node[below]{$\delta_i$}(0,0.8)node[above]{$B_k$};

\draw[red,thick,->-=.5,>=stealth] (0,0)\nn to (2,0);

	\end{tikzpicture}
\end{minipage}

\caption{$(\SS,\MM,\TT)$ to $(\SS',\MM',\TT')$}\label{f:Some changes in self-folded triangle}
\end{figure}

\section{Multiset of  arcs determined by their dimension vectors}\label{s:Multiset of untagged arcs determined by their dimension vectors}

\subsection{Intersection numbers}
The definition of intersection number for arcs was introduced by \cite[Definition 8.4]{FST08}

\begin{df}\cite{FST08} Let $\alpha$ and $\beta$ be two  arcs in $(\SS,\MM)$. The intersection number $\dInt(\alpha|\beta)$ is defined as follows:
$$\dInt(\alpha|\beta)= \Int^A(\alpha|\beta)+\Int^B(\alpha|\beta)+\Int^C(\alpha|\beta),$$
where
\begin{itemize}
\item $\Int^A(\alpha|\beta)$ is the number of intersection points of $\alpha$ and $\beta$ in $(\SS,\MM)$;
\item $\Int^B(\alpha|\beta)=0$ unless $\alpha$ is a loop based at a marked point $a$, in which case $\Int^B(\alpha|\beta)$ is computed
as follows: assume that $\beta$ intersects $\alpha$ (in the interior of $\SS\setminus \MM$) at the points $b_1,\cdots,b_m$ (numbered along $\beta$ in this order); then $\Int^B(\alpha|\beta)$ is the negative of the number of segments $\gamma_i=[b_i, b_{i+1}]\subset \beta$ such that $\gamma_i$ together with the segments $[a, b_i]\subset \alpha$ and $[a, b_{i+1}]\subset \alpha$ form three sides of a contractible triangle disjoint from the punctures;
\item $\Int^C(\alpha|\beta)=0$ unless $\alpha=\beta$, in which case $\Int^C(\alpha|\beta)=-1$;
\end{itemize}
\end{df}

\subsection{Intersection vectors}
For an ideal triangulation $\TT$ of $(\SS,\MM)$,  we define the {\it intersection vector} of an arc $\alpha$ with respect to $\TT$ as: $$\dIntv_{\TT}(\alpha):=(\dInt(\ba |\alpha))_{\ba\in \TT}\in \mathbb{Z}^n,$$
where $n$ is the cardinality of $\fT$.
Let $\cm$ be a multiset of pairwise compatible arcs, define $\dIntt(\ba|\cm)=\Sigma_{\alpha\in\cm}\Int(\ba|\alpha)$ and $$\dIntv_{\TT}(\cm):=\Sigma_{\alpha\in\cm}\dIntv_{\TT}(\alpha).$$
We call $\dIntv_{\TT}(\cm)$ the {\it intersection vector} of $\cm$ with respect to $\TT$.
Note that $\dIntv_{\TT}(\cm)$ is also equal to $(\dInt(\ba|\cm))_{\ba\in\TT},$ \ie,  $$\dIntv_{\TT}(\cm)=(\dInt(\ba|\cm))_{\ba\in\TT}.$$

\subsection{Notation}
Let $\Delta$  be a triangle, $P$ a vertex of $\Delta$ and $\gamma$ an edge of $\Delta$. Additionally, let $\angle \eta$ be an angle within $\Delta$. For a finite multiset $\cm$ consisted of pairwise compatible arcs, we need the following notations:
\begin{itemize}
    \item $\Arcseg^1_{\Delta}(\cM)$: the multiset of all  the irreducible arc segments  of $\cM$ lying in $\Delta$;
    \item  $\Arcseg^1_{\Delta}(\cm, \angle \eta)$: the multiset of all the irreducible arc segments in  $\Arcseg^1_{\Delta}(\cM)$ that cross $\angle \eta$;
    \item $\Arcseg^1_{\Delta}(\cm, P, \gamma)$:  the multiset  of all  the irreducible arc segments   in  $\Arcseg^1_{\Delta}(\cM)$, having one  endpoint at $P$ and the other  in interior of $\gamma$;
    \item $|S|$: the count of $S$ for any multiset $S$.
\end{itemize}

\subsection{Local equivalence via the same intersection vectors}
In this subsection, let $\TT$ be a strong admissible triangulation of $(\SS,\MM)$,  $\cm$ and $\cn$ be two multisets consisting of pairwise compatible arcs, such that $\dIntv_{\TT}(\cm)=\dIntv_{\TT}(\cn)$. Moreover, we assume that 
 neither $\cm$ nor $\cn$ contain any arcs of $\TT$. Consequently, for any arc $\gamma$ of $\cm$ or  $\cn$ and  $\ba\in\TT$, $\dInt(\ba|\gamma)=\Int^A(\ba|\gamma)+\Int^B(\ba|\gamma)$.

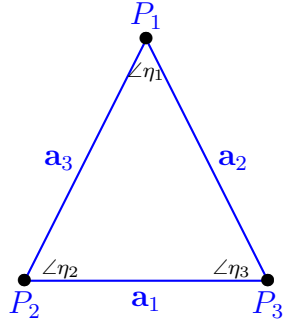
\begin{figure}
	
		\begin{tikzpicture}[xscale=0.8,yscale=0.8]
			
			\draw[blue,thick] (-2,-2)\nn to (2,-2)\nn;
			
		\draw[blue,thick] (-2,-2)\nn to (0,2)\nn;
			\draw[blue,thick] (2,-2)\nn to (0,2)\nn;

		\draw[blue](-2,-2)node[below]{$P_2$}(2,-2)node[below]{$P_3$}(0,2)node[above]{$P_1$};	
		\draw[blue](0,-2)node[below]{$\ba_1$}(1,0)node[right]{$\ba_2$}(-1,0)node[left]{$\ba_3$};			
	\draw[black](-1.9,-1.8)node[right]{{\tiny$\angle\eta_2$}}(1.9,-1.8)node[left]{{\tiny$\angle\eta_3$}}(0,1.8)node[below]{{\tiny$\angle\eta_1$}};		
			
		\end{tikzpicture}
	
\caption{Proper triangle}\label{f:tri-gons}	
\end{figure}

\begin{lem}\label{l:general case}
Let $\Delta$ be a proper triangle, as depicted in Figure~\ref{f:tri-gons}. If $\Int^A(\ba_i|\cm) = \Int^A(\ba_i|\cn)$ holds for any  $1 \leq i \leq 3$, then  $\Arcseg^1_{\Delta}(\cM)=\Arcseg^1_{\Delta}(\cN)$. 
\end{lem}

\begin{proof}
Let us denote $a_i=|\Arcseg^1_{\Delta}(\cm,P_i,\ba_i)|$, $b_i=|\Arcseg^1_{\Delta}(\cm,\angle \eta_i)|$, $a_i'=|\Arcseg^1_{\Delta}(\cn,P_i,\ba_i)|$ and $b_i'=|\Arcseg^1_{\Delta}(\cn,\angle \eta_i)|$ for each $1\leq i\leq 3$.
It  suffices to prove that $a_i=a_i'$ and $b_i=b_i'$ for each $1\leq i\leq 3$.

We begin by demonstrating that $a_i\neq 0$ if and only if $a_i'\neq 0$ for each $1\leq i\leq 3$.
 Without loss of generality, let us assume  $a_1\neq 0$ but $a_1'= 0$. 
 By   $a_1\neq 0$, we have $\Int^A(\ba_1|\cm)=a_1+b_2+b_3$, $\Int^A(\ba_2|\cm)=b_3$, $\Int^A(\ba_3|\cm)=b_2$ (cf. the first picture of Figure~\ref{f:Cases of no edge of delta is loop}).
 Hence, we have $$\Int^A(\ba_1|\cm)-\Int^A(\ba_2|\cm)-\Int^A(\ba_3|\cm)=a_1>0.$$
By   $a_1'= 0$,
one can deduce that $$\Int^A(\ba_1|\cn)-\Int^A(\ba_2|\cn)-\Int^A(\ba_3|\cn)\leq 0$$ (cf. the last two  pictures of  Figure \ref{f:Cases of no edge of delta is loop}). This leads to a contradiction, since $\Int^A(\ba_i|\cm) = \Int^A(\ba_i|\cn)$  for every  $1 \leq i \leq 3$. Then by  discussing the following two cases,  we can finish our proof.

\noindent\textbf{Case} 1:  $a_i\neq 0$ for some $1\leq i\leq 3$. Without loss of generality, assume  $a_1\neq 0$, hence $a_1'\neq 0$.  We then have 
\[\Int^A(\ba_1|\cm)=a_1+b_2+b_3,\Int^A(\ba_2|\cm)=b_3, \Int^A(\ba_3|\cm)=b_2
\]
as well  as 
\[\Int^A(\ba_1|\cn)=a_1'+b_2'+b_3',\Int^A(\ba_2|\cn)=b_3',\Int^A(\ba_3|\cn)=b_2'.
\] 
By  $\Int^A(\ba_i|\cm)=\Int^A(\ba_i|\cn)$ for each $\ba_i\in\TT $, 
we have
\[
\left\{
\begin{array}{ccc}
a_1+b_2+b_3&=&a_1'+b_2'+b_3',\\
                b_2&=&b_2',\\
               b_3&=&b_3'.
\end{array}
\right.
\]
It is easy to get that  $a_i=a_i'$, $b_i=b_i'$ for each $1\leq i\leq 3$.

\noindent\textbf{Case} 2: $a_i= 0$ for each $1\leq i\leq 3$. So $a_i'= 0$ for each $1\leq i\leq 3$ as well.
 Then 
 $\Int^A(\ba_1|\cm)=b_2+b_3$, $\Int^A(\ba_2|\cn)=b_1+b_3$ and $\Int^A(\ba_3|\cn)=b_1+b_2$ (cf. the second  pictures of  Figure \ref{f:Cases of no edge of delta is loop}). Similarly, 
  $\Int^A(\ba_1|\cn)=b_2'+b_3'$, $\Int^A(\ba_2|\cn)=b_1'+b_3'$, and $\Int^A(\ba_3|\cn)=b_1'+b_2'$.  
Then \[
\left\{
\begin{array}{ccc}
b_2+b_3&=&b_2'+b_3',\\
b_1+b_3&=&b_1'+b_3',\\
 b_1+b_2&=& b_1'+b_2'.
\end{array}
\right.
\]
It is easy to get that  $a_i=a_i'$, $b_i=b_i'$ for each $1\leq i\leq 3$.
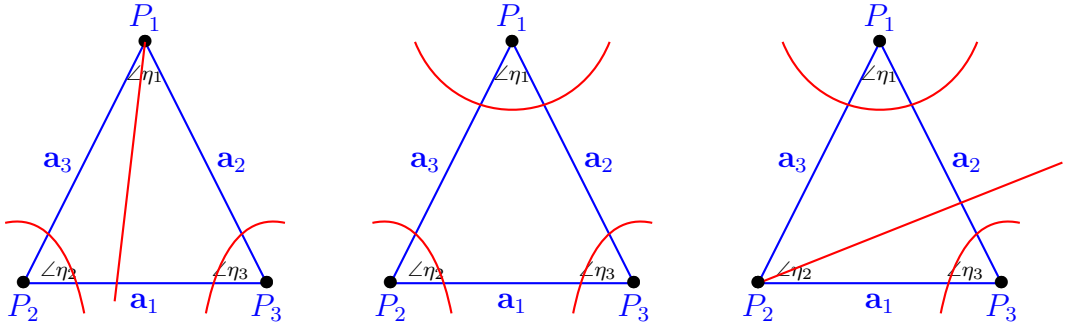
\begin{figure}[h]

	\begin{minipage}{0.3\linewidth} 
		\begin{tikzpicture}[xscale=0.8,yscale=0.8]
			
			\draw[blue,thick] (-2,-2)\nn to (2,-2)\nn;
			
		\draw[blue,thick] (-2,-2)\nn to (0,2)\nn;
			\draw[blue,thick] (2,-2)\nn to (0,2)\nn;

		\draw[blue](-2,-2)node[below]{$P_2$}(2,-2)node[below]{$P_3$}(0,2)node[above]{$P_1$};	
		\draw[black](-1.9,-1.8)node[right]{{\tiny$\angle\eta_2$}}(1.9,-1.8)node[left]{{\tiny$\angle\eta_3$}}(0,1.8)node[below]{{\tiny$\angle\eta_1$}};

		\draw[blue](0,-2)node[below]{$\ba_1$}(1,0)node[right]{$\ba_2$}(-1,0)node[left]{$\ba_3$};			
		
	\draw[red,thick] (0,2) to (-0.5,-2.3);
\draw[thick,color=red] (-2.3,-1) .. controls (-1.4,-0.8) and (-1.1,-2) .. (-1,-2.5);	
\draw[thick,color=red] (2.3,-1) .. controls (1.4,-0.8) and (1.1,-2) .. (1,-2.5);	

		\end{tikzpicture}
		\end{minipage}
\begin{minipage}{0.3\linewidth}
		\begin{tikzpicture}[xscale=0.8,yscale=0.8]
			
			\draw[blue,thick] (-2,-2)\nn to (2,-2)\nn;
			
		\draw[blue,thick] (-2,-2)\nn to (0,2)\nn;
			\draw[blue,thick] (2,-2)\nn to (0,2)\nn;

		\draw[blue](-2,-2)node[below]{$P_2$}(2,-2)node[below]{$P_3$}(0,2)node[above]{$P_1$};	
	\draw[black](-1.9,-1.8)node[right]{{\tiny$\angle\eta_2$}}(1.9,-1.8)node[left]{{\tiny$\angle\eta_3$}}(0,1.8)node[below]{{\tiny$\angle\eta_1$}};

		\draw[blue](0,-2)node[below]{$\ba_1$}(1,0)node[right]{$\ba_2$}(-1,0)node[left]{$\ba_3$};

\draw[thick,color=red] (-2.3,-1) .. controls (-1.4,-0.8) and (-1.1,-2) .. (-1,-2.5);	

\draw[thick,color=red] (2.3,-1) .. controls (1.4,-0.8) and (1.1,-2) .. (1,-2.5);	
	
\draw[thick,color=red] (-1.6,2) .. controls (-1,0.5) and (1,0.5) .. (1.6,2);	
		\end{tikzpicture}
		
		\end{minipage}
	\begin{minipage}{0.3\linewidth}
		\begin{tikzpicture}[xscale=0.8,yscale=0.8]
			
			\draw[blue,thick] (-2,-2)\nn to (2,-2)\nn;
			
		\draw[blue,thick] (-2,-2)\nn to (0,2)\nn;
			\draw[blue,thick] (2,-2)\nn to (0,2)\nn;

		\draw[blue](-2,-2)node[below]{$P_2$}(2,-2)node[below]{$P_3$}(0,2)node[above]{$P_1$};	
	\draw[black](-1.9,-1.8)node[right]{{\tiny$\angle\eta_2$}}(1.9,-1.8)node[left]{{\tiny$\angle\eta_3$}}(0,1.8)node[below]{{\tiny$\angle\eta_1$}};

		\draw[blue](0,-2)node[below]{$\ba_1$}(1,0)node[right]{$\ba_2$}(-1,0)node[left]{$\ba_3$};			
		
\draw[red,thick] (-2,-2)\nn to (3,0);

\draw[thick,color=red] (2.3,-1) .. controls (1.4,-0.8) and (1.1,-2) .. (1,-2.5);	
	
\draw[thick,color=red] (-1.6,2) .. controls (-1,0.5) and (1,0.5) .. (1.6,2);	
		\end{tikzpicture}
		
		\end{minipage}
	\caption{$\Delta$ is  proper triangle}\label{f:Cases of no edge of delta is loop}
	\end{figure}
	\end{proof}

\begin{lem}\label{l:tri-gons}
Let $\Delta$ be a $3$-vertices triangle (cf. Figure~\ref{f:Basic triangle under strong admissible triangulation}), 
then  $\Arcseg^1_{\Delta}(\cM)=\Arcseg^1_{\Delta}(\cN)$. 
\end{lem}
\begin{proof}
Suppose that the vertices and edges of $\Delta$ are depicted in Figure~\ref{f:tri-gons}. 
Given that no edge of $\Delta$ is a loop, for any arc $\gamma\in\cm$ or $\gamma\in\cn$, we have $\Int^B(\ba_i|\gamma)=0$ for $1\leq i\leq 3$, hence $\dInt(\ba_i|\gamma)=\Int^A(\ba_i|\gamma)$. Therefore, by $\dInt(\ba_i|\cm)=\dInt(\ba_i|\cn)$, it follows that $\Int^A(\ba_i|\cm)=\Int^A(\ba_i|\cn)$. Consequently,
by Lemma~\ref{l:general case}, we can get the  result.
\end{proof}

\begin{figure}[h]
\begin{minipage}{0.4\linewidth}
\begin{tikzpicture}[xscale=0.8,yscale=0.8]
			\draw[blue,thick] (0,-2)\nn to (0,0)\nn;
			
\draw[thick,color=blue] (0,-2) .. controls (-2,-1.5) and (-2,2) .. (0,2);	
\draw[thick,color=blue] (0,-2) .. controls (2,-1.5) and (2,2) .. (0,2);		
\draw[thick,blue](0,-0.6)node[right]{$\ba_P$}(0,2)node[above]{$\bl_P$}(0,-2)node[below]{$Q$}(0,0)node[left]{$P$};
\draw[thick,black](-0.4,-1.8)node[above]{{\tiny$\angle \eta_1$}}(0.5,-1.8)node[above]{{\tiny$\angle \eta_2$}};
		
	\end{tikzpicture}
\end{minipage}
\begin{minipage}{0.4\linewidth}
\begin{tikzpicture}[xscale=0.8,yscale=0.8]
			\draw[blue,thick] (0,-2)\nn to (0,0)\nn;
			
\draw[thick,color=blue] (0,-2) .. controls (-2,-1.5) and (-2,2) .. (0,2);	
\draw[thick,color=blue] (0,-2) .. controls (2,-1.5) and (2,2) .. (0,2);		
\draw[thick,blue](0,-0.6)node[right]{$\ba_P$}(0,2)node[above]{$\bl_P$}(0,-2)node[below]{$Q$}(0,0)node[left]{$P$};
\draw[thick,black](-0.4,-1.8)node[above]{{\tiny$\angle \eta_1$}}(0.5,-1.8)node[above]{{\tiny$\angle \eta_2$}};
\draw[thick,color=red] (0,0) .. controls (0.5,0.5) and (0.8,1.2) .. (1,2);	
\draw[thick,color=red] (-1,-2) .. controls (-0.5,-1.5) and (0.5,-1.5) .. (1,-2);

	\end{tikzpicture}
\end{minipage}

\caption{Self-folded triangle}\label{f:Self-folded triangle}
\end{figure}
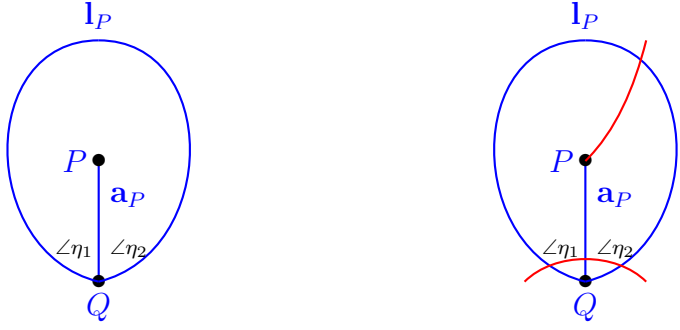

\begin{lem}\label{l:self-folded triangle}
Let $\Delta$ be a self-folded triangle, then $\Arcseg^1_{\Delta}(\cM)=\Arcseg^1_{\Delta}(\cN)$.
\end{lem}
\begin{proof}
Suppose the vertices, edges and angles of $\Delta$ are depicted 
 in Figure~\ref{f:Self-folded triangle}, where $Q$ is a boundary marked point, $P$ is the puncture, $\ba_P$ is the folded side, $\bl_P$ is the remaining side.  
It suffices to prove that   $|\Arcseg^1_{\Delta}(\cm, \angle \eta_i)|=|\Arcseg^1_{\Delta}(\cn, \angle \eta_i)|$ for each $1\leq i\leq 2$  and  $|\Arcseg^1_{\Delta}(\cm, P, \bl_P)|= |\Arcseg^1_{\Delta}(\cn, P, \bl_P)|$.
Note that $|\Arcseg^1_{\Delta}(\cm, \angle \eta_1)|=|\Arcseg^1_{\Delta}(\cm, \angle \eta_2)|$, $|\Arcseg^1_{\Delta}(\cn, \angle \eta_1)|=|\Arcseg^1_{\Delta}(\cn, \angle \eta_2)|$, 
 \[
\left\{
\begin{array}{ccc}
\dInt(\ba_P|\cm)&=&|\Arcseg^1_{\Delta}(\cm,\angle\eta_1)|,\\
\dInt(\bl_P|\cm)&=&|\Arcseg^1_{\Delta}(\cm,\angle\eta_1)|+|\Arcseg^1_{\Delta}(\cm,P, l_P)|,
\end{array}
\right.
\]
(cf. Figure~\ref{f:Self-folded triangle}) and 
 \[
\left\{
\begin{array}{ccc}
\dInt(\ba_P|\cn)&=&|\Arcseg^1_{\Delta}(\cn,\angle\eta_1)|,\\
\dInt(
\bl_P|\cn)&=&|\Arcseg^1_{\Delta}(\cn,\angle\eta_1)|+|\Arcseg^1_{\Delta}(\cm,P, 
\bl_P)|.
\end{array}
\right.
\]
By the condition that
 $\dInt(\ba_P|\cm)=\dInt(\ba_P|\cn)$ and $\dInt(\bl_P|\cm)=\dInt(\bl_P|\cn)$,  one can obtain that $|\Arcseg^1_{\Delta}(\cm, \angle \eta_i)|=|\Arcseg^1_{\Delta}(\cn, \angle \eta_i)|$ for each $1\leq i\leq 2$  and  $|\Arcseg^1_{\Delta}(\cm, P, \bl_P)|= |\Arcseg^1_{\Delta}(\cn, P, \bl_P)|$. 
\end{proof}

\begin{figure}[htpb]

		\begin{tikzpicture}

			\draw[blue,thick] (0,-2)\nn to (0,-2)\nn;
				\draw[blue,thick] (0,1)\nn to (0,1)\nn;

\draw[thick,color=blue] (0,-2) .. controls (-1,-1.5) and (-1,0) .. (0,0);	
\draw[thick,color=blue] (0,-2) .. controls (1,-1.5) and (1,0) .. (0,0);	


\draw[thick,color=blue] (0,-2) .. controls (-3,-1.5) and (-1,1) .. (0,1);	
\draw[thick,color=blue] (0,-2) .. controls (3,-1.5) and (1,1) .. (0,1);	
		
\draw[blue](0,-2)node[below]{{\small$Q$}}(0,1.2)node[above]{{\small$P_1$}}(0,-1)node[above]{{\small$P$}};	
\draw[blue](-1.3,0)node[left]{{\small$\ba_3$}}(1.3,0)node[right]{{\small$\ba_2$}}(0,0)node[below]{{\small$\ba_1$}};			
\draw[black](0,0.8)node{{\tiny$\angle\eta_1$}}(-0.9,-1.5)node{{\tiny$\angle\eta_2$}}	(0.9,-1.5)node{{\tiny$\angle\eta_3$}}(-0.9,0)node{{\small$\Delta$}};	
\draw[thick,black](-0.3,-1.6)node[above]{{\tiny$\angle \eta_4$}}(0.3,-1.6)node[above]{{\tiny$\angle \eta_5$}};
	\draw[blue,thick] (0,-2)\nn to (0,-1)\nn;

		\end{tikzpicture}

		\caption{ $\Delta$ is  a 2-vertices triangle}\label{f: edge of delta is a  loop}
	\end{figure}
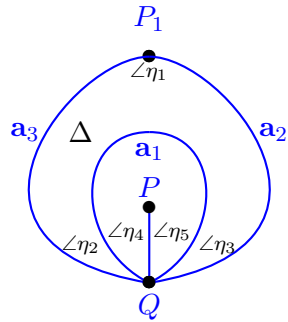

\begin{lem}\label{l: triangle of type I with one edge is a remaining side}
Let $\Delta$ be a $2$-vertices triangle (cf. Figure~\ref{f:Basic triangle under strong admissible triangulation}), then $\Arcseg^1_{\Delta}(\cM)=\Arcseg^1_{\Delta}(\cN)$.
\end{lem}

\begin{proof} Suppose the vertices, edges and angles of $\Delta$ are
depicted in Figure~\ref{f: edge of delta is a  loop},  where $\ba_1$ is just the remaining side of a self-folded triangle, $P$ is a puncture, $P_1,Q$ are boundary points such that $P_1\neq Q$.
Then $\Int^B(\ba_1|\cm)=-|\Arcseg^1_{\Delta}(\cm,\angle \eta_4)|$ and $\Int^B(\ba_1|\cn)=-|\Arcseg^1_{\Delta}(\cn,\angle \eta_4)|$.
According to  Lemma~\ref{l:self-folded triangle}, $$|\Arcseg^1_{\Delta}(\cm,\angle \eta_4)|=|\Arcseg^1_{\Delta}(\cn,\angle \eta_4)|,$$ which implies that $\Int^B(\ba_1|\cm)=\Int^B(\ba_1|\cn)$.
Since $\ba_1$  is a loop, we have
$$\Int^A(\ba_1|\cm)=\dInt(\ba_1|\cm)-\Int^B(\ba_1|\cm)=\dInt(\ba_1|\cn)-\Int^B(\ba_1|\cn)= \Int^A(\ba_1|\cn).$$
Because $\ba_2$ and $\ba_3$ are not loops,
 we have $$\Int^A(\ba_i|\cm)=\dInt(\ba_i|\cm)=\dInt(\ba_i|\cn)=\Int^A(\ba_i|\cn)$$
 for $i=2,3$.
 We conclude the result by Lemma~\ref{l:general case}.
   \end{proof}

\subsection{Multisets of untagged arcs determined by their intersection vectors }
For an arc $\gamma$ and a positive integer $m$, we denote by $\gamma^{(m)}$ the multiset that contains $m$ copies of $\gamma$. Moreover, we denote by $\gamma^{(0)}=\emptyset$.
\begin{lem}\label{l: T has no element in cm and cn} Let $\cm,\cn$ be two multisets consisting of pairwise compatible arcs such that $\dIntv_{\TT}(\cm)=\dIntv_{\TT}(\cn)$. 
Let $\ba \in \TT$ be an arc.  Assume that $\ba^{(m)} \subseteq \cm, \ba^{(m+1)}\not\subseteq\cm$ and $\ba^{(n)} \subseteq \cn, \ba^{(n+1)}\not\subseteq\cn$ for some non-negative integers $m$ and $n$, then $m = n$.

\end{lem}
\begin{proof} 
Without loss of generality, we may assume that $m>0$. We are going to show $m=n$.
Given that all arcs  in $\cm$ are compatible,
   we have $\dInt(\ba|\gamma)=0$ for any $\gamma\in \cm$ with $\gamma\neq \ba$. On the other hand, $\dInt(\ba|\ba)=-1$.
It follows that $\dInt(\ba|\cm)=-m$. 
Since $\dIntv_{\TT}(\cm)=\dIntv_{\TT}(\cn)$, we conclude that $\dInt(\ba|\cn)=-m$.  We claim that $n>0$. Otherwise $n=0$. In other words, $\ba\not\in \cn$. According to the definition of intersection number,  $\dInt(\ba|\cn)\geq 0$, which contradicts to $\dInt(\ba|\cn)=-m<0$. Hence $n>0$.
A similar discussion yields $\dInt(\ba|\cn)=-n$ and hence $m=n$.


 
\end{proof}
Now we can state the main result of this section.
\begin{thm}\label{t:case of s has three vertex}
Let $(\SS,\MM)$ be a marked surface with at least three marked points on the boundary,  $\TT$  a strong admissible ideal triangulation of $(\SS,\MM)$.  Let $\cm,\cn$ be two multisets consisting of pairwise compatible arcs respectively such that $\dIntv_{\TT}(\cm)=\dIntv_{\TT}(\cn)$, then $\cm=\cn$.
\end{thm}
\begin{proof}
By Lemma~\ref{l: T has no element in cm and cn}, we can assume neither $\cm$ nor $\cn$ contains arcs from $\TT$. 
Consider a triangle $\Delta$, because  $\TT$ is a strong admissible ideal triangulation of $(\SS,\MM)$, then $\Delta$ is either a self-folded triangle, a 3-vertices triangle  or 2-vertices triangle. According to Lemmas \ref{l:tri-gons}, \ref{l:self-folded triangle}, and \ref{l: triangle of type I with one edge is a remaining side},
we have $\Arcseg^1(\cM)=\Arcseg^1(\cN)$.
Consequently, it follows that $\cm=\cn$ by Theorem \ref{t:local to entire}. 
\end{proof}

\section{Multiset of tagged arcs determined by their intersection vectors}\label{s:Multiset of tagged arcs determined by their intersection vectors}
\subsection{Tagged arcs}
Let $(\SS,\MM)$ be a marked surface.
 Each arc $\gamma$ in $(\SS, \MM)$ has two ends obtained by arbitrarily cutting $\gamma$ into three pieces, and then throwing out the middle one. A {\it tagged arc} is an arc in which each end has been tagged in one of two ways, plain or notched, so that the following conditions are satisfied:
 \begin{itemize}
\item[($\mathbf{T1}$)] the arc does not cut out a once-punctured monogon; 
\item[($\mathbf{T2}$)] an endpoint lying on the boundary is tagged plain; 
\item[($\mathbf{T3}$)] both ends of a loop are tagged in the same way.
\end{itemize}

For a tagged arc $\gamma$,
 \begin{itemize}
\item $\gamma$ is called a {\it plain arc} if its both ends are tagged plain;
\item $\gamma$ is called a {\it 1-notched arc} if  one end of $\gamma$ is tagged plain and the other end is tagged notched;
\item $\gamma$ is called a {\it 2-notched arc} if its both ends are tagged notched.
\end{itemize}
In the figures, we represent tags as follows:
\begin{figure}[h]
\begin{minipage}[t]{0.3\linewidth} 
\begin{tikzpicture}[xscale=0.6,yscale=0.6]
\draw[blue,thick] (-4,0) to (-2,0)\nn;
\draw[blue,thick](2,0) to (4,0)\nn;
\node at (-5,0){plain};
\node at (0.5,0){notched};
\node at (3.5,0)[rotate=90]{$\bowtie$};
\end{tikzpicture}
\end{minipage}
\end{figure}

For a tagged arc $\alpha$ in $(\SS,\MM)$, we denote by 
 $\overline{\alpha}$ the untagged version of $\alpha$. For tagged arcs $\alpha$ and $\beta$ such
that $\overline{\alpha}=\overline{\beta}$, if exactly one of them is a 1-notched arc, then the pair $(\alpha,\beta)$ is called {\it a pair of conjugate arcs}.  By $(\mathbf{T3})$, if  $(\alpha,\beta)$ is a  pair of conjugate arcs, neither $\alpha$ nor $\beta$ is a loop.
 For a puncture $P$, if there is a pair of conjugate arcs which have different tags at the ends incident to  $P$, we can denote by $(\gamma_P^-,\gamma_P^{\bowtie})$ such conjugate arcs, where $\gamma_P^-$ is tagged plain at the end incident to $P$ while $\gamma_P^{\bowtie}$ is tagged notched at the end incident to $P$. 

\subsection{Tagged triangulation}
Two tagged arcs $\alpha,\beta$ are called {\it compatible} if and only if the following conditions are satisfied:
\begin{itemize}
\item[($\mathbf{C1}$)] $\overline{\alpha}$ and $\overline{\beta}$ are compatible;
\item[($\mathbf{C2}$)] if $\overline{\alpha}\neq \overline{\beta}$ and $\alpha$ and $\beta$ share an endpoint $P$, then the ends of $\alpha$ and $\beta$ connected to $P$ must be tagged in the same way;
\item[($\mathbf{C3}$)] if $\overline{\alpha}=\overline{\beta}$, then at least one end of $\alpha$ must be tagged in the same way as the corresponding end of $\beta$.
\end{itemize}

\noindent For each pair of conjugate  arcs $(\gamma_P^-,\gamma_P^{\bowtie})$ at some puncture $P$, one can infer that $\gamma_P^-$ and $\gamma_P^{\bowtie}$ are compatible.

A maximal  set of pairwise compatible tagged arcs in $(\SS,\MM)$ is
called a {\it tagged triangulation} of $(\SS,\MM)$.

\begin{figure}[htpb]		
 \begin{minipage}{0.35\linewidth}
		\begin{tikzpicture}[xscale=0.7,yscale=0.7]
			
			\draw[blue,thick] (-2,-2)\nn to (2,-2)\nn;
			
		\draw[blue,thick] (-2,-2)\nn to (0,2)\nn;
			\draw[blue,thick] (2,-2)\nn to (0,2)\nn;

		
		\draw[blue](0,-3)node[below]{Triangle piece};
		\end{tikzpicture}
		
	\end{minipage}		
 \begin{minipage}{0.3\linewidth}
		\begin{tikzpicture}[xscale=0.7,yscale=0.7]

		\draw[blue,thick] (0,2)\nn to (0,2)\nn;
			\draw[blue,thick] (0,-2)\nn to (0,-2)\nn;
				\draw[blue,thick] (0,0)\nn to (0,0)\nn;
\draw[thick,color=blue] (0,-2) .. controls (-2,-1) and (-2,1) .. (0,2);	
\draw[thick,color=blue] (0,-2) .. controls (2,-1) and (2,1) .. (0,2);	

\draw[thick,color=blue] (0,-2) .. controls (-0.3,-1.4) and (-0.3,-0.6) .. (0,0);	
\draw[thick,color=blue] (0,-2) .. controls (0.3,-1.4) and (0.3,-0.6) .. (0,0);

		
\draw[blue](0.17,-0.4)node[rotate=20]{$\bowtie$};			
\draw[blue](0,-3)node[below]{1-puncture piece};
		\end{tikzpicture}
\end{minipage}		
 \begin{minipage}{0.3\linewidth}
		\begin{tikzpicture}[xscale=0.7,yscale=0.7]

			\draw[blue,thick] (0,-2)\nn to (0,-2)\nn;
			\draw[blue,thick] (-1,0)\nn to (-1,0)\nn;
				\draw[blue,thick] (1,0)\nn to (1,0)\nn;
\draw[thick,color=blue] (0,-2) .. controls (-2,-1.5) and (-2,2) .. (0,2);	
\draw[thick,color=blue] (0,-2) .. controls (2,-1.5) and (2,2) .. (0,2);

\draw[thick,color=blue] (0,-2) .. controls (-0.2,-0.8) and (-0.65,0) .. (-1,0);	
\draw[thick,color=blue] (0,-2) .. controls (-0.6,-1.3) and (-1,-0.8) .. (-1,0);	

\draw[blue](-0.63,-0.23)node[rotate=40]{$\bowtie$};

\draw[thick,color=blue] (0,-2) .. controls (0.2,-0.8) and (0.65,0) .. (1,0);	
\draw[thick,color=blue] (0,-2) .. controls (0.6,-1.3) and (1,-0.8) .. (1,0);	

\draw[blue](0.63,-0.23)node[rotate=-40]{$\bowtie$};

\draw[blue](0,-3)node[below]{2-puncture piece};		
		\end{tikzpicture}
  \end{minipage}
		\caption{Basic tiles}\label{f:basic tiles}
	\end{figure}
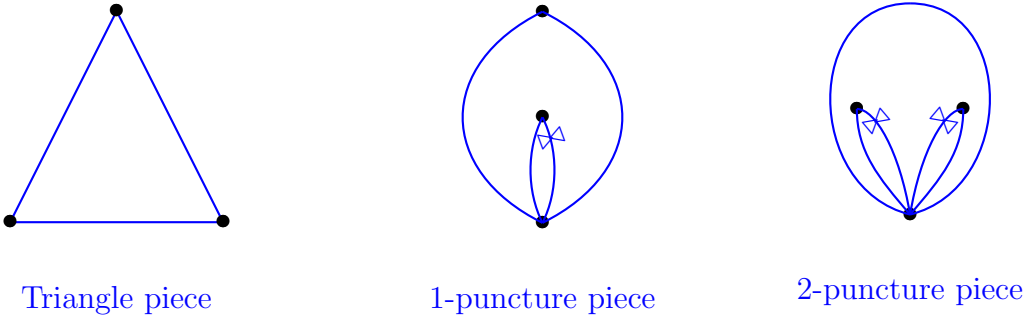

Let $\TT$ be a tagged triangulation. Given that $\SS$ is assumed to have a boundary, the basic tiles are of three types: the first is a triangle with three distinct edges, the second is a digon containing a puncture and a pair of conjugate arcs within, referred to as a {\it $1$-puncture piece}, and the third is a monogon with two punctures and two pairs of conjugate arcs inside, known as a {\it $2$-puncture piece}.
 (cf. Figure~\ref{f:basic tiles}).
\subsection{(Strong) admissible tagged triangulation} \label{s:(Strong) admissible tagged triangulation}
A tagged triangulation $\TT$ is called {\it admissible} if each puncture of $(\SS,\MM)$ is connected to a boundary marked point by a pair of conjugate arcs in $\TT$. Because all tagged arcs in $\TT$ are compatible with each other, for any puncture, beside the pair of conjugate arcs of $\TT$ linked to it,  there are no additional arcs in $\TT$ that connect to this puncture.
So for an admissible  triangulation $\TT$, the puncture appears either in the interior of $1$-puncture piece or in the interior of $2$-puncture piece (cf. Figure~\ref{f:basic tiles}) and for each puncture $P$, there is a unique pair of conjugate arcs in $\TT$ incident to $P$, where we denote this pair by $(\ba_P^-,\ba_P^{\bowtie})$ in the following.

An admissible tagged triangulation $T$ is termed {\it strong} if $\TT$ does not include any loops. Consequently, in a strong admissible tagged triangulation, the basic tile is either a triangle piece with three distinct boundary marked points serving as vertexes or a $1$-puncture piece with two distinct boundary marked points as vertexes.

\subsection{Signature of a partial tagged triangulation} Let $C$ be a multiset of pairwise compatible tagged arcs. For a puncture $P$, let
us temporarily define $S_C(P)\subseteq\{\mathbf{plain, notched}\}$ as the collection of tags of tagged arcs  within  $C$ appearing at the end incident to $P$.
It should be observed that if $S_C(P) = \{\mathbf{plain,notched}\}$, then there exists a pair of conjugate arcs $(\gamma^{-}_P,\gamma^{\bowtie}_P)$ in $C$ incident to $P$. Any other tagged arcs in $C$ that are incident to $P$ are either a copy of $\gamma^{-}_P$ or a copy of $\gamma^{\bowtie}_P$.


\subsection{Tagged triangulation and ideal triangulation}
Let $C$ be a set of pairwise compatible tagged arcs. 
We define another set $C^{\circ}$, which consists of plain arcs derived from $C$ through the following process:
\begin{itemize}
\item  For each puncture $P$ where $S_C(P) = \{\mathbf{notched}\}$, we replace all (notched) tags at $P$ by plain ones.
\item For each puncture $P$ with $S_C(P) = \{\mathbf{plain,notched}\}$, replace the arc $\beta$ notched at $P$  by the corresponding remaining side $l_P$ enclosing $P$ and $\beta$; that is, $l_P$ is a loop based at the endpoint $Q\neq P$ of $\beta$,  is obtained
by closely wrapping around $\beta$.
\end{itemize}
Note that if $S_C(P) = \{\mathbf{plain,notched}\}$ at some puncture $P$, then there is  a pair of conjugate arcs $(\gamma_P^{-},\gamma_P^{\bowtie})$ incident to $P$, and no other tagged arcs in $C$ connected to $P$ since any two tagged arcs in $C$ are pairwise compatible  and $C$ is a set. So $\beta=\gamma_P^{\bowtie}$ and $\gamma_P^-$ is closely enclosed by  $l_P$. 
 We denote this replacement as $(\gamma^{\bowtie}_P)^{\circ} = l_P$.
It is straightforward to check that arcs in $C^{\circ}$ are pairwise compatible.

The following is a direct consequence of definitions.

\begin{lem}
   Let $(\SS,\MM)$ be a marked surface with a (strong) admissible tagged triangulation $\TT$. Then  $\TT^\circ$ is  a  (strong) admissible ideal triangulation of $(\SS,\MM)$.
\end{lem}


Let $\TT$ be an admissible tagged triangulation and $\ba\in \TT$, we denote \[
\ba^{\circ}=\left\{
\begin{array}{ccc}
\bl_P&\text{if } \ba=\ba_P^{\bowtie} \text{ for some puncture } P;\\
\ba&\text{else}.
\end{array}
\right.
\]

\subsection{Intersection numbers of tagged arcs}
Following  \cite{FST08}, we generalize the intersection number of plain arcs to tagged arcs.
\begin{df}\cite{FST08} Let $\alpha$ and $\beta$ be two tagged arcs in $(\SS,\MM)$. The intersection number $\Int(\alpha|\beta)$ is defined as follows:
$$\dIntt(\alpha|\beta)= \dInt(\overline{\alpha}|\overline{\beta})+\Int^D(\alpha|\beta),$$
where
\begin{itemize}

\item $\Int^D(\alpha|\beta)$ is the number of ends of $\beta$ that are incident to an endpoint of $\alpha$ and, at that endpoint, carry a tag different from the one $\alpha$ does. 
\end{itemize}
\end{df}
\noindent For convenience of notation, we also denote by   $\Int^A(\alpha|\beta):=\Int^A(\overline{\alpha}|\overline{\beta})$, $\Int^B(\alpha|\beta):=\Int^B(\overline{\alpha}|\overline{\beta})$,  $\Int^C(\alpha|\beta):=\Int^C(\overline{\alpha}|\overline{\beta})$ and $\dInt(\alpha|\beta):=\dInt(\overline{\alpha}|\overline{\beta})$. 
So for any tagged arcs $\alpha$ and $\beta$, we also have $$\dInt(\alpha|\beta)= \Int^A(\alpha|\beta)+\Int^B(\alpha|\beta)+\Int^C(\alpha|\beta).$$
 It is clear that if  $\Int^D(\alpha|\beta)=0$, then $\Int(\alpha|\beta)=\dInt(\alpha|\beta)$.
For a pair of conjugate arcs $(\ba_P^-,\ba_P^{\bowtie})$ and any tagged arc $\beta$, we have 
$\Int^{*}(\ba_P^-|\beta)=\Int^{*}(\ba_P^{\bowtie}|\beta)$ for $*\in\{A,B,C, \circ\}$.

\subsection{Intersection vectors of tagged arcs}
For a tagged triangulation $\TT$ of $(\SS,\MM)$,  we define the {\it intersection vector} of a tagged arc $\alpha$ with respect to $\TT$ as follows: $$\Intv_{\TT}(\alpha):=(\Int(\ba |\alpha))_{\ba\in \TT}\in \Z_{\ge 0}^n,$$
where $n=|\TT|$.
Let $\cm$ be a multiset of pairwise compatible tagged arcs, we set  $\dIntt(\ba|\cm)=\Sigma_{\alpha\in\cm}\Int(\ba|\alpha)$ and define the {\it intersection vector of $\cm$} with respect to $\TT$ as  $$\Intv_{\TT}(\cm):=\Sigma_{\alpha\in\cm}\Intv_{\TT}(\alpha).$$
Note that $\Intv_{\TT}(\cm)$ is also equal to $(\Int(\ba|\cm))_{\ba\in\TT},$ \ie,  $$\Intv_{\TT}(\cm)=(\Int(\ba|\cm))_{\ba\in\TT}.$$

\begin{lem}\label{l:relation of l_p and a_p}
 Let $\TT$ be an admissible tagged triangulation of $(\SS,\MM)$. 
 Let $P$ be a puncture, $(\ba_P^-,\ba_P^{\bowtie})$ the corresponding pair of conjugate arcs in $\TT$. Recall that $\bl_P=(\ba_P^{\bowtie})^\circ$ is the loop wrapping $\ba_P^-$.
  Let $\alpha$ be a tagged arc of $(\SS,\MM)$,  then  
\begin{itemize}
\item [(1)] $\Int^B(\bl_P|\alpha) = -\Int^A(\ba_P^{-}|\alpha).$
    \item[(2)] If neither endpoint of $\alpha$ is $P$, then $$\Int^D(\ba_P^{\bowtie}|\alpha) = \Int^D(\ba_P^{-}|\alpha) = 0 \ \text{and}\ \Int^A(\bl_P|\alpha) = 2\Int^A(\ba_P^{-}|\alpha).$$
    \item[(3)] If only one endpoint of $\alpha$ is $P$, then $\Int^D(\ba_P^{\bowtie}|\alpha) + \Int^D(\ba_P^{-}|\alpha) = 1$. Moreover, if $\alpha\neq \ba_P^{-}$ and $\alpha\neq \ba_P^{\bowtie}$,  we have   $$\Int^A(\bl_P|\alpha) = 2\Int^A(\ba_P^{-}|\alpha) + 1.$$
\item[(4)] If both endpoints of $\alpha$ are $P$, then $$\Int^D(\ba_P^{\bowtie}|\alpha) + \Int^D(\ba_P^{-}|\alpha) = 2   \ \text{and}\ \Int^A(\bl_P|\alpha) = 2\Int^A(\ba_P^{-}|\alpha) + 2.$$
\item[(5)] For any $\alpha\notin\TT$, we have
\begin{align*}
  \Int^A(\bl_P|\alpha)+\Int^B(\bl_P|\alpha)&=\Int^A(\ba_P^{\bowtie}|\alpha)+\Int^D(\ba_P^{\bowtie}|\alpha)+\Int^D(\ba_P^{-}|\alpha) \\
&=\Int^A(\ba_P^{-}|\alpha)+\Int^D(\ba_P^{\bowtie}|\alpha)+\Int^D(\ba_P^{-}|\alpha). 
\end{align*}
\end{itemize}

\end{lem}

\begin{proof} 
For $(1)$, because $\bl_P$ is a one-puntured loop and $\bl_P$  closely wraps around $\ba_P^-$, by the definitions of $\Int^B(\bl_P|\alpha)$ and $\Int^A(\ba_P^{-}|\alpha)$, one can get $\Int^B(\bl_P|\alpha) = -\Int^A(\ba_P^{-}|\alpha)$ directly.

Since $\TT$ is admissible, the  other endpoint of $\ba_P^{-}$ and $\ba_P^{\bowtie}$ is not a puncture.  Hence when compute $\Int^D(\ba_P^{\bowtie}|\alpha)$ and $ \Int^D(\ba_P^{-}|\alpha)$, we only need to check the  end of $\alpha$ that is incident to $P$. 
If  neither endpoint of $\alpha$ is $P$, it is clear that $\Int^D(\ba_P^{\bowtie}|\alpha) = \Int^D(\ba_P^{-}|\alpha) = 0$. Else if only one endpoint of $\alpha$ is $P$, so $\alpha$ carries different tag  either with $\ba_P^{\bowtie}$ or with $ \ba_P^{-}$, hence  $\Int^D(\ba_P^{\bowtie}|\alpha) + \Int^D(\ba_P^{-}|\alpha) = 1$. If both endpoints of $\alpha$ are $P$, so $\alpha$ is a loop, thus both ends of $\alpha$ are tagged in the same way, which is also different  either with $\ba_P^{\bowtie}$ or with $ \ba_P^{-}$, then $\Int^D(\ba_P^{\bowtie}|\alpha) + \Int^D(\ba_P^{-}|\alpha) = 2$. This proves the first parts in $(2)$,$(3)$ and $(4)$.

We say  the field enclosed by $\bl_P$, \ie the field where $P$ lies in, is the inside of $\bl_P$. 
Because $\ba_P^-$ is closely wrapped by $\bl_P$, so we can regard $\bl_P$ as $(\ba_P^-)^{(2)}$ when $\ba_P^-$ passes through  $\alpha$. In other words,   we can regard $\bl_P$ as $(\ba_P^-)^{(2)}$ when $\alpha$ passes through $\ba_P^-$  from outside  of $\bl_P$ to outside  of $\bl_P$.
During this process, one intersection of $\ba_P^-$ with $\alpha$ is accompanied with two intersections of $\bl_P$  with $\alpha$. Note that, all intersections of   $\ba_P^-$ and $\alpha$ are created in this way. So we only need to consider the intersections of $\bl_P$  with $\alpha$ which are not created in this way. 

Because $\bl_P$ is the  loop that closely wraps around $\ba_P^-$, so if neither endpoint of $\alpha$ is $P$, then no  endpoints of $\alpha$  are inside $\bl_P$, hence, no more intersections of $\bl_P$  with $\alpha$ are created,
then $\Int^A(\bl_P|\alpha) = 2\Int^A(\ba_P^{-}|\alpha)$. Else if  only one endpoint of $\alpha$ is $P$, $\alpha\neq \ba_P^{-}$ and $\alpha\neq \ba_P^{\bowtie}$, then $\alpha$ must pass through $\bl_P$. Thus,
one more intersection of  $\bl_P$  with $\alpha$ is created when $\alpha$ passes   from inside  of $\bl_P$ to outside  of $\bl_P$ with $P$ as the starting point, 
so 
$\Int^A(\bl_P|\alpha) = 2\Int^A(\ba_P^{-}|\alpha) + 1$.  If both endpoints of $\alpha$ are $P$, two more intersections of  $\bl_P$  with $\alpha$ are created
when  $\alpha$  passes  from inside  of $\bl_P$ to outside  of $\bl_P$ at the beginning and  $\alpha$  passes from outside  of $\bl_P$ to inside  of $\bl_P$ at the ending. It follows that 
$\Int^A(\bl_P|\alpha) = 2\Int^A(\ba_P^{-}|\alpha) + 2$. This completes the proof of $(2),(3)$ and $(4)$.

By noticing that $\Int^A(\ba_P^-|\alpha)=\Int^A(\ba_P^{\bowtie}|\alpha)$, we obtain $(5)$ by $(1)$--$(4)$.

\end{proof}

\begin{lem}\label{l:relation of a_p and l_Q}
  Let $\TT$ be an admissible triangulation of $(\SS,\MM)$.  Let $(\gamma_Q^-,\gamma_Q^{\bowtie})$ be a pair of conjugate arcs  at the puncture $Q$, $l_Q=(\gamma_Q^{\bowtie})^{\circ}$ the loop closely wrapping $\gamma_Q^-$. The following statements hold.
  \begin{itemize}
      \item[(1)] 
      For each $\ba\in\TT$,  $$\Int^B(\ba|l_Q)=2\Int^B(\ba|\gamma_Q^{-}).$$
 \item[(2)] Let $P$ be a puncture, $(\ba_P^-,\ba_P^{\bowtie})$ the corresponding pair of conjugate arcs in $\TT$, $\bl_P=(\ba_P^{\bowtie})^\circ$ the loop wrapping $\ba_P^-$. \begin{itemize} 
  \item[(2.1)] If $P\neq Q$, then
$$\Int^B(\bl_P|l_Q)=2\Int^B(\bl_P|\gamma_Q^-).$$
\item[(2.2)] If $P=Q$, then $$\Int^B(\bl_P|l_Q)=2\Int^B(\bl_P|\gamma_Q^-)-1.$$
\item[(2.3)] We have $$\Int^A(\bl_P|l_Q)=2\Int^A(\bl_P|\gamma_Q^-)=\Int^A(\bl_P|\gamma_Q^-)+\Int^A(\bl_P|\gamma_Q^{\bowtie}).$$ 
     \end{itemize}

      \item[(3)] Let $(\ba_Q^-,\ba_Q^{\bowtie})$ be the pair of conjugate  arcs in $\TT$ corresponding to $Q$, $\ba\in\TT$. 
  \begin{itemize}
      \item [(3.1)]  If $\ba\neq \ba_Q^-$ and $\ba\neq \ba_Q^{\bowtie}$,
 then $$\Int^A(\ba|l_Q)=2\Int^A(\ba|\gamma_Q^{-})=\Int^A(\ba|\gamma_Q^{-})+\Int^A(\ba|\gamma_Q^{\bowtie}).$$
 \item[(3.2)] If $\ba=\ba_Q^-$ or $\ba=\ba_Q^{\bowtie}$, assume moreover that $(\gamma_Q^-,\gamma_Q^{\bowtie})\neq (\ba_Q^-,\ba_Q^{\bowtie})$, then 
$$\Int^A(\ba|l_Q)=2\Int^A(\ba|\gamma_Q^{-})+1.$$ 
  \end{itemize} 
 \end{itemize}
  
\end{lem}
\begin{proof}


Let  $\bl$ be a loop  based at a boundary marked point $Z$.   Because $l_Q$ is the loop closely wrapping $\gamma_Q^-$,  when  $\gamma_Q^-$ passes through $\bl$, we can regard $l_Q$ as $(\gamma_Q^-)^{(2)}$.  In particular, when $\gamma_Q^-$ passes through from outside of $\bl$ to inside of $\bl$, then to outside of $\bl$, one segment of $\gamma_Q^-$ within $\bl$ is created, two same  segments of $l_Q$ within $\bl$ are accompanied created. Note that all segments of $\gamma_Q^-$  within $\bl$  are created in this way. So we only need to consider the segments of $l_Q$ within $\bl$ which are not created in this way. 
If $Q$ is not enclosed by $\bl$, then one can get that no more segment of  $l_Q$  within $\bl$   is created. It follows that $\Int^B(\bl|l_Q)=2\Int^B(\bl|\gamma_Q^{-})$.
If $Q$ is enclosed by $\bl$,  when $\gamma_Q^-$ first passes from inside of $\bl$ to outside of $\bl$ with $Q$ as the starting point, one more segment  of  $l_Q$ within $\bl$  is created.  Denote this segment by $[A, B]$, where 
$A,B$ are the intersections of $\bl$ with $l_Q$. 

For the part $(1)$, if $\ba$ is not a loop, then $\Int^B(\ba|l_Q)=2\Int^B(\ba|\gamma_Q^{-})=0$ by definition. 
Now assume that $\ba=\bl\in\TT$ is a loop. If $Q$ is not  enclosed by $\ba$, so $\Int^B(\ba|l_Q)=2\Int^B(\ba|\gamma_Q^{-})$ by above discussion. If $Q$ is   enclosed by $\ba$, so  $\ba$ is not a once-punctured loop by $(\mathbf{T1})$, then the triangle formed by the segments $[Z,A]\subset \ba$, $[Z,B]\subset \ba$ and $[A,B]\subset l_Q$ is not a contractible triangle. Hence, we have 
$\Int^B(\ba|l_Q)=2\Int^B(\ba|\gamma_Q^{-})$. This finishes the proof of $(1)$.

For $(2.1)$, let $\bl=\bl_P$, because  $P\neq Q$,  $Q$ is not enclosed by $\bl_P$. It follows that $$\Int^B(\bl_P|l_Q)=2\Int^B(\bl_P|\gamma_Q^-).$$

For $(2.2)$, let $\bl=\bl_P$, because $P=Q$,  $Q$ is enclosed by $\bl_P$. Then the triangle formed by the segments $[Z,A]\subset \bl_P$, $[Z,B]\subset \bl_P$ and $[A,B]\subset l_Q$ is  a contractible triangle, so $$\Int^B(\bl_P|l_Q)=2\Int^B(\bl_P|\gamma_Q^-)-1.$$

For $(2.3)$, because $Q$ is a puncture and the base point of $\bl_P$ is not a puncture, we deduce that the base point of $\bl_P$ is not $Q$. By using same proof as $(2)$ of Lemma~\ref{l:relation of l_p and a_p}, we conclude that $\Int^A(l_Q|\bl_P)=2\Int^A(\gamma_Q^-|\bl_P)$. Hence, $$\Int^A(\bl_P|l_Q)=\Int^A(l_Q|\bl_P)=2\Int^A(\gamma_Q^-|\bl_P)=2\Int^A(\bl_P|\gamma_Q^-)=\Int^A(\bl_P|\gamma_Q^-)+\Int^A(\bl_P|\gamma_Q^{\bowtie}).$$ 

For $(3)$,  
 if $\ba\neq \ba_Q^-$ and $\ba\neq \ba_Q^{\bowtie}$, then $Q$ is not an endpoint of $\ba$. By $(2)$ of Lemma~\ref{l:relation of l_p and a_p}, we have $$\Int^A(\ba|l_Q)=2\Int^A(\ba|\gamma_Q^{-}).$$ 
Else if $\ba=\ba_Q^-$ or $\ba=\ba_Q^{\bowtie}$,  only one endpoint of $\ba$ is $Q$, because $(\gamma_Q^-,\gamma_Q^{\bowtie})\neq (\ba_Q^-,\ba_Q^{\bowtie})$, by using same proof as $(3)$ of Lemma~\ref{l:relation of l_p and a_p}, we have
$$\Int^A(\ba|l_Q)=2\Int^A(\ba|\gamma_Q^{-})+1.$$
\end{proof}

\subsection{ Modification of a tagged triangulation}  Consider a multiset $C$ of pairwise compatible tagged arcs, we define a function $\sigma_C(P)$ as the difference between the number of ends of arcs in $C$ incident to   $P$ that are tagged plain and the number of ends of arcs in $C$ incident to   $P$ that  are tagged notched, \ie, let $(\ba_P^-,\ba_P^{\bowtie})$ be a pair of conjugate arcs, then $$\sigma_C(P)=\Int^D(\ba^{\bowtie}|C)-\Int^D(\ba^{-}_P|C).$$
 A positive value  of $\sigma_C(P)$  indicate more plain ends than notched ends  of arcs in $C$ incident to   $P$, while a negative value indicate the opposite.

Given an admissible tagged triangulation $\TT$ of $(\SS,\MM)$ and considering a multiset $C$ of pairwise compatible tagged arcs, we proceed to construct a new triangulation $\TT_C$ from $\TT$ by carrying out the following steps:
 for each puncture $P$ with $\sigma_{C}(P) < 0$, replace  $\ba_{P}^{-}$ and  $\ba_{P}^{\bowtie}$  in $\TT$ with each other.
Given this construction, there exists a bijection $\phi_C$ from $\TT$ to $\TT_C$ such that:
\begin{itemize}
    \item $\phi_C(\ba_{P}^{-}) = \ba_{P}^{\bowtie}$ and $\phi_C(\ba_{P}^{\bowtie}) = \ba_{P}^{-}$ when $\sigma_{C}(P) < 0$.
    \item $\phi_C(\ba) = \ba$ for all other tagged arcs in $\TT$.
\end{itemize}
It is evident that $\TT_C$ is also an admissible tagged triangulation. Moreover, $\TT_C$ is strong admissible if $\TT$ is strong admissible.

\begin{lem}\label{l:arc between tag and untag1}
Let $\TT$ be  an admissible tagged triangulation of $(\SS,\MM)$ and 
  a multiset  $\CM$ of pairwise compatible tagged arcs such that $\CM\cap \TT=\emptyset$ and $S_{\CM}(P)\neq \{\mathbf{plain,notched}\}$ for each puncture $P$.
   Let $\alpha$  be an  tagged arc in $\CM$,  
   we have 
   $$\dInttv_{\TT}(\alpha)=\dIntv_{\TT_{\CM}^{\circ}}(\alpha)=\dIntv_{\TT_{\CM}^{\circ}}(\overline{\alpha}).$$ 
 \end{lem}

\begin{proof}  Recall that $\overline{\alpha}$ is the untagged version of $\alpha$, it is clear that $\dIntv_{\TT_{\CM}^{\circ}}(\alpha)=\dIntv_{\TT_{\CM}^{\circ}}(\overline{\alpha})$. Then
it suffices to prove that  for each $\ba\in \TT$, $\Int(\ba|\alpha)=\dInt((\phi_{\CM}(\ba))^\circ|\alpha).$
Since ${\CM}\cap\TT=\emptyset$, it follows that $\Int^{C}(\ba|\alpha) = 0$ for any arc $\ba$ in $\TT$ and  $\alpha$ in ${\CM}$. Therefore, 
$$\dInt(\ba|\alpha) = \Int^A(\ba|\alpha) + \Int^B(\ba|\alpha),\ \ \Int(\ba|\alpha) = \Int^A(\ba|\alpha) + \Int^B(\ba|\alpha) + \Int^D(\ba|\alpha).$$
Because $\TT$ is admissible,  either neither endpoint of $\ba$ is a puncture or only one endpoint of $\ba$ is a puncture. Therefore, we only need to consider the following two cases.

\noindent$\mathbf{Case}$ 1: Neither endpoint of $\ba$ is a puncture. Then $(\phi_{\CM}(\ba))^\circ = \ba^\circ = \ba$ and $\Int^D(\ba|\alpha) = 0$. Thus, 
$$\dInt((\phi_{\CM}(\ba))^\circ|\alpha) = \dInt(\ba|\alpha) = \Int(\ba|\alpha).$$

\noindent$\mathbf{Case}$ 2: One endpoint of $\ba$ is a puncture, say $P$. Because $\TT$ is admissible, there is a unique pair $(\ba_P^-,\ba_P^{\bowtie})$  of conjugate arcs in $\TT$ incident to $P$. Hence $\ba = \ba_P^{-}$ or $\ba_P^{\bowtie}$. It follows that  $\ba$ is not  a loop. Consequently,
$\Int^B(\ba|\alpha) = 0$.  Hence,
\begin{align}\label{gongshi:intersection of a and alpha in case 2}
 \dInt(\ba|\alpha) = \Int^A(\ba|\alpha),\ \ \Int(\ba|\alpha) = \Int^A(\ba|\alpha) + \Int^D(\ba|\alpha)
\end{align}
 in this case. Denote by $\bl_P=(\ba_P^{\bowtie})^{\circ}$,
note that $\alpha$ cannot cut out an once punctured loop by $(\mathbf{T1})$, so $\overline{\alpha}\neq \bl_P$, which implies $\Int^C(\bl_P|\alpha)=0,$ then $$\dInt(\bl_P|\alpha)=\Int^A(\bl_P|\alpha)+\Int^B(\bl_P|\alpha)+\Int^C(\bl_P|\alpha)=\Int^A(\bl_P|\alpha)+\Int^B(\bl_P|\alpha).$$
Let $Z$ be the base point of $\bl_P$, so $Z$ is also the endpoint of $\ba_P^{-}$ and $\ba_P^{\bowtie}$ and $Z\neq P$.
\noindent$\mathbf{Subcase}$ 2.1: 
Suppose first $\sigma_P({\CM})\ge 0$, which implies $\phi_{\CM}(\ba)=\ba$.  Since $S_{\CM}(P)\neq \{\mathbf{plain,notched}\}$, we have $S_{\CM}(P)=\{\mathbf{plain}\}$ or $S_{\CM}(P)=\emptyset$. As a consequence, each arc  in ${\CM}$ with $P$ as an endpoint is tagged plain at the end incident to $P$. Because $Z$ is not a puncture,  we  have 
  $\Int^D(\ba_P^{-}|\alpha)=0$. 
Therefore,   
$$\dInt((\phi_{\CM}(\ba_P^{-}))^\circ|\alpha)=\dInt(\ba_P^{-}|\alpha)=\Int(\ba_P^{-}|\alpha)$$
and 
\begin{align*}
\dInt((\phi_{\CM}(\ba_P^{\bowtie}))^\circ|\alpha)
&=\dInt(\bl_P|\alpha)\\
&=\Int^A(\bl_P|\alpha)+\Int^B(\bl_P|\alpha)\\
&=\Int^A(\ba_P^{\bowtie}|\alpha)+\Int^D(\ba_P^{-}|\alpha)+\Int^D(\ba_P^{\bowtie}|\alpha)& &( \text{by (5) of Lemma}~\ref{l:relation of l_p and a_p})\\
&=\Int^A(\ba_P^{\bowtie}|\alpha)+\Int^D(\ba_P^{\bowtie}|\alpha)\\
&=\Int(\ba_P^{\bowtie}|\alpha)& &(\text{by} (\ref{gongshi:intersection of a and alpha in case 2})).
\end{align*}

\noindent$\mathbf{Subcase}$ 2.2: Now suppose  $\sigma_P({\CM})< 0$, then $\phi(\ba_P^{\bowtie})=\ba_P^-$, $\phi(\ba_P^{-})=\ba_P^{\bowtie}$ and $S_{\CM}(P)=\{\mathbf{notched}\}$.
If $P$ is an endpoint of $\alpha$,  since all the arcs in ${\CM}$ are pairwise compatible and $S_{\CM}(P)=\{\mathbf{notched}\}$, then $\alpha$ is tagged notched at the end incident to $P$. It follows that $\Int^D(\ba_P^{\bowtie}|\alpha)=0$. If $P$ is not an endpoint of $\alpha$, because $Z$ is not a puncture, then $\ba_P^{\bowtie}$ and $\alpha$ do not have a common endpoint which is a puncture. As a consequence, we also have $\Int^D(\ba_P^{\bowtie}|\alpha)=0$ in this situation.
 Then by (\ref{gongshi:intersection of a and alpha in case 2}), we have
 \begin{eqnarray*}
\dInt((\phi_{\CM}(\ba_P^{\bowtie}))^\circ|\alpha)
=\dInt(\ba_P^{-}|\alpha)=\Int^A(\ba_P^{-}|\alpha)
=\Int(\ba_P^{\bowtie}|\alpha).
\end{eqnarray*}
On the other hand,
\begin{align*}\dInt((\phi_{\CM}(\ba_P^{-}))^\circ|\alpha)&=\dInt(\bl_P|\alpha)\\
&=\Int^A(\bl_P|\alpha)+\Int^B(\bl_P|\alpha)\\
&=\Int^A(\ba_P^{-}|\alpha)+\Int^D(\ba_P^{-}|\alpha)+\Int^D(\ba_P^{\bowtie}|\alpha)& &( \text{by (5) of Lemma}~\ref{l:relation of l_p and a_p})\\
&=\Int^A(\ba_P^{-}|\alpha)+\Int^D(\ba_P^{-}|\alpha)\\
&=\Int(\ba_P^{-}|\alpha) & &(\text{by } (\ref{gongshi:intersection of a and alpha in case 2})).
\end{align*}

\end{proof}

\begin{lem}\label{l:arc between tag and untag2}
Let $\TT$ be an admissible tagged triangulation of $(\SS,\MM)$ and ${\CM}$
  a multiset  of pairwise compatible tagged arcs such that ${\CM}\cap \TT=\emptyset$.  
  Let $(\gamma_{Q}^{-},\gamma_{Q}^{\bowtie})$  be a pair  of  conjugate  arcs
 in ${\CM}$, denote by $l_{Q}=(\gamma_{Q}^{\bowtie})^\circ$  the loop closely wrapping around $\gamma_{Q}^{-}$,
then under the bijection $\phi_{\CM}$, we have $$\dInttv_{\TT}(\gamma_{Q}^{-})+\dInttv_{\TT}(\gamma_{Q}^{\bowtie})=\dIntv_{\TT_{{\CM}}^{\circ}}(l_{Q}).$$

 \end{lem}

\begin{proof}
It suffices to prove that  for each $\ba\in \TT$, $$\Int(\ba|\gamma_{Q}^-)+\Int(\ba|\gamma_{Q}^{\bowtie})=\dInt((\phi_{\CM}(\ba))^\circ|l_{Q}).$$
Since ${\CM}\cap\TT=\emptyset$, it follows that $\Int^C(\ba|\gamma) = 0$ for any arc $\ba\in\TT$ and $\gamma\in {\CM}$. Therefore, we can express  $\Int(\ba|\gamma)$ as:
\begin{align}\label{gongshi:intersection of a and gamma in general case}
    \Int(\ba|\gamma) = \Int^A(\ba|\gamma) + \Int^B(\ba|\gamma) + \Int^D(\ba|\gamma).
\end{align}
 Let ${P}$ be a puncture. Since $\TT$ is admissible,  there is a unique pair of conjugate arcs $(\ba_{P}^-,\ba_{P}^{\bowtie})$ incident to $P$, denote by $\bl_{P}=(\ba_{P}^{\bowtie})^\circ$ the loop wrapping $\ba_{P}^-$. 
By  ${\CM}\cap\TT=\emptyset$, we conclude that $\bl_P\neq l_Q$. Hence, 
$\Int^C(\bl_P|l_Q) = 0.$ It follows that 
\begin{align}\label{gongshi:intersection of lp and lq in general case}
\dInt(\bl_P|l_Q)=\Int^A(\bl_P|l_Q)+\Int^B(\bl_P|l_Q)+\Int^C(\bl_P|l_Q)=\Int^A(\bl_P|l_Q)+\Int^B(\bl_P|l_Q).   
\end{align}
Note that each arc in $\TT$ is not a once-punctured loop. We have $\Int^C(\ba|l_Q) = 0$  for each $\ba\in\TT$.
Therefore
\begin{align}\label{gongshi:intersection of a and lq in general case}
\dInt(\ba|l_Q)=\Int^A(\ba|l_Q)+\Int^B(\ba|l_Q)+\Int^C(\ba|l_Q)=\Int^A(\ba|l_Q)+\Int^B(\ba|l_Q).  
\end{align}
In particular,
\begin{align}\label{gongshi:intersection of ap- and lq in general case}
\dInt(\ba_P^{-}|l_Q)=\Int^A(\ba_P^{-}|l_Q).  
\end{align}

\noindent$\mathbf{Case}$ 1: No endpoint of $\ba$ is a puncture. Then $\Int^D(\ba|\gamma_{Q}^{-})=\Int^D(\ba|\gamma_{Q}^{\bowtie})=0$. By (\ref{gongshi:intersection of a and gamma in general case}), we have      $$\Int(\ba|\gamma_{Q}^{-})+\Int(\ba|\gamma_{Q}^{\bowtie})=\Int^A(\ba|\gamma_{Q}^{-})+\Int^A(\ba|\gamma_{Q}^{\bowtie})+\Int^B(\ba|\gamma_{Q}^{-})+\Int^B(\ba|\gamma_{Q}^{\bowtie}).$$
It follows that
\begin{align*} \dInt((\phi_{\CM}(\ba))^\circ|l_{Q})
&=\dInt(\ba|l_{Q})\\
&=\Int^A(\ba|l_{Q})+\Int^B(\ba|l_{Q}) \\
&=2\Int^A(\ba|\gamma_{Q}^{-})+2\Int^B(\ba|\gamma_{Q}^{-})\ \ \ \ \ \ \ \ \ \ (\text{by Lemma}~\ref{l:relation of a_p and l_Q})\\
&=\Int^A(\ba|\gamma_{Q}^{-})+\Int^A(\ba|\gamma_{Q}^{\bowtie})+\Int^B(\ba|\gamma_{Q}^{-})+\Int^B(\ba|\gamma_{Q}^{\bowtie})\\
&=\Int(\ba|\gamma_{Q}^{-})+\Int(\ba|\gamma_{Q}^{\bowtie}).
\end{align*}

\noindent$\mathbf{Case}$ 2: One endpoint of $\ba$ is a puncture $P$ and $P\neq Q$. Since $\TT$ is admissible, then
 $\ba=\ba_{P}^{-}$ or $\ba_{P}^{\bowtie}$. In particular, $\ba$ is not a loop and hence $$\Int^B(\ba|l_{Q})=\Int^B(\ba|\gamma_{Q}^{\bowtie})=\Int^B(\ba|\gamma_{Q}^{-})=0.$$ By (\ref{gongshi:intersection of a and gamma in general case}),  we have
 \begin{align}\label{gongshi: intersection of a and gammaQ in case 2}
   \Int(\ba|\gamma_{Q}^{-})+\Int(\ba|\gamma_{Q}^{\bowtie})=\Int^A(\ba|\gamma_{Q}^{\bowtie})+\Int^A(\ba|\gamma_{Q}^{-})+\Int^D(\ba|\gamma_{Q}^{-})+\Int^D(\ba|\gamma_{Q}^{\bowtie}).  
 \end{align}
 On the other hand, we also have
\begin{align}\label{gongshi: intersection of a and lQ in case 2}
    \dInt(\ba|l_Q)=\Int^A(\ba|l_Q)
\end{align}
by (\ref{gongshi:intersection of a and lq in general case}).
Furthermore, 
\begin{align*} 
&\dInt(\bl_{P}|l_{Q})\\
&=\Int^A(\bl_{P}|l_{Q})+\Int^B(\bl_{P}|l_{Q})\quad \quad \quad \quad \quad\quad \quad \quad \quad \quad(\text{by (\ref{gongshi:intersection of lp and lq in general case})})\\
 &=2\Int^A(\bl_{P}|\gamma_{Q}^{-})+2\Int^B(\bl_{P}|\gamma_{Q}^{-})\quad \quad \quad \quad \quad\quad \quad \quad (\text{by (2.3) and (3.1) of Lemma \ref{l:relation of a_p and l_Q}})\\
  &=\Int^A(\bl_{P}|\gamma_{Q}^{-})+\Int^A(\bl_{P}|\gamma_{Q}^{\bowtie})+\Int^B(\bl_{P}|\gamma_{Q}^{-})+\Int^B(\bl_{P}|\gamma_{Q}^{\bowtie})\\
&=\Int^A(\ba_{P}^{\bowtie}|\gamma_{Q}^{-})+\Int^D(\ba_{P}^{\bowtie}|\gamma_{Q}^{-})+\Int^D(\ba_{P}^{-}|\gamma_{Q}^{-})+ \Int^A(\ba_{P}^{\bowtie}|\gamma_{Q}^{\bowtie})+\Int^D(\ba_{P}^{\bowtie}|\gamma_{Q}^{\bowtie})+\Int^D(\ba_{P}^{-}|\gamma_{Q}^{\bowtie}),
\end{align*}
where  the last equality follows from (5) of Lemma~\ref{l:relation of l_p and a_p}.

Let $Z$ (resp. $O$) be the base point of $\bl_P$ (resp. $l_Q$). Hence $Z$ (resp. $O$) is also the endpoint of $\ba_P^{-}$ and $\ba_P^{\bowtie}$ (resp. $\gamma_Q^{-}$ and $\gamma_Q^{\bowtie}$) and $Z\neq P$ (resp. $O\neq Q$).

\noindent$\mathbf{Subcase}$ 2.1:  $P\neq O$.
 Since $Z$ is not a puncture and $P\neq Q$, then  $\Int^D(\ba|\gamma_{Q}^{-})=\Int^D(\ba|\gamma_{Q}^{\bowtie})=0$. Hence, by (\ref{gongshi: intersection of a and gammaQ in case 2}), we have
\begin{align*} 
\Int(\ba|\gamma_{Q}^{-})+\Int(\ba|\gamma_{Q}^{\bowtie})=\Int^A(\ba|\gamma_{Q}^{-})+\Int^A(\ba|\gamma_{Q}^{\bowtie}).
\end{align*}
Moreover, we have 
\begin{align*} 
\dInt(\ba_{P}^{-}|l_{Q})&=\Int^A(\ba_{P}^{-}|l_{Q})&& (\text{by (\ref{gongshi: intersection of a and lQ in case 2})})\\
&=\Int^A(\ba_{P}^{-}|\gamma_{Q}^{-})+\Int^A(\ba_{P}^{-}|\gamma_{Q}^{\bowtie})&&(\text{by (3.1) of Lemma \ref{l:relation of a_p and l_Q}})
\end{align*}
and 
 \begin{align*} 
&\dInt(\bl_{P}|l_{Q})\\
=&\Int^A(\ba_{P}^{\bowtie}|\gamma_{Q}^{-})+\Int^D(\ba_{P}^{\bowtie}|\gamma_{Q}^{-})+\Int^D(\ba_{P}^{-}|\gamma_{Q}^{-})+ \Int^A(\ba_{P}^{\bowtie}|\gamma_{P}^{\bowtie})+\Int^D(\ba_{P}^{\bowtie}|\gamma_{Q}^{\bowtie})+\Int^D(\ba_{P}^{-}|\gamma_{Q}^{\bowtie})\\
=&\Int^A(\ba_{P}^{\bowtie}|\gamma_{Q}^{-})+\Int^A(\ba_{P}^{\bowtie}|\gamma_{Q}^{\bowtie}).
\end{align*}
Since $\Int^A(\ba_P^{\bowtie}|\gamma_{Q}^{-})=\Int^A(\ba_P^{\bowtie}|\gamma_{Q}^{\bowtie})=\Int^A(\ba_P^{-}|\gamma_{Q}^{-})=\Int^A(\ba_P^{-}|\gamma_{Q}^{\bowtie})$,  we obtain 
$$\dInt(\ba_P^{-}|l_{Q})=\dInt(\bl_{P}|l_{Q})=\Int(\ba_P^{\bowtie}|\gamma_{Q}^{-})+\Int(\ba_P^{\bowtie}|\gamma_{Q}^{\bowtie})=\Int(\ba_P^{-}|\gamma_{Q}^{-})+\Int(\ba_P^{-}|\gamma_{Q}^{\bowtie}).$$
Hence, if $\sigma_{P}({\CM})< 0$, then
\begin{align*}
\dInt((\phi_{\CM}(\ba_P^{\bowtie}))^\circ|l_{Q})=\dInt(\ba_P^{-}|l_{Q})
=\Int(\ba_P^{\bowtie}|\gamma_{Q}^{-})+\Int(\ba_P^{\bowtie}|\gamma_{Q}^{\bowtie}),
\end{align*}
\begin{align*}
\dInt((\phi_{\CM}(\ba_P^{-}))^\circ|l_{Q})=\dInt(\bl_{P}|l_{Q})
=\Int(\ba_P^{-}|\gamma_{Q}^{-})+\Int(\ba_P^{-}|\gamma_{Q}^{\bowtie}).
\end{align*}
Similar discussion yields the result for $\sigma_{P}({\CM})\ge 0$.

\noindent$\mathbf{Subcase}$ 2.2: $P=O$.
First, assume that
$\sigma_{P}({\CM})\ge 0$, then  $\gamma_Q^-$ and $\gamma_Q^{\bowtie}$ are tagged plain at the end incident to $P$. It follows that
\begin{align}\label{gongshi: D-intersection of a and gammaq in subcase 2.2.1}
    \Int^D(\ba_P^-|\gamma_{Q}^{-})=\Int^D(\ba_P^{-}|\gamma_{Q}^{\bowtie})=0.
\end{align}
According to (\ref{gongshi: intersection of a and gammaQ in case 2}), we have 
\begin{align}\label{gongshi: intersection of ap- and gammaq in subcase 2.2.1}
    \Int(\ba_{P}^{-}|\gamma_{Q}^{-})+\Int(\ba_{P}^{-}|\gamma_{Q}^{\bowtie})=\Int^A(\ba_{P}^{-}|\gamma_{Q}^{-})+\Int^A(\ba_{P}^{-}|\gamma_{Q}^{\bowtie}).
\end{align}
Consequently,
\begin{align*} 
\dInt((\phi_{\CM}(\ba_{P}^{-}))^\circ|l_{Q})&=\dInt(\ba_{P}^{-}|l_{Q})\\
&=\Int^A(\ba_{P}^{-}|l_{Q})&& (\text{by }~(\ref{gongshi: intersection of a and lQ in case 2}))\\
&=2\Int^A(\ba_{P}^{-}|\gamma_{Q}^{-})&& (\text{by (3.1) of Lemma}~\ref{l:relation of a_p and l_Q})\\
&=\Int^A(\ba_{P}^{-}|\gamma_{Q}^{-})+\Int^A(\ba_{P}^{-}|\gamma_{Q}^{\bowtie})\\
&=\Int(\ba_{P}^{-}|\gamma_{Q}^{-})+\Int(\ba_{P}^{-}|\gamma_{Q}^{\bowtie})&&(\text{by } (\ref{gongshi: intersection of ap- and gammaq in subcase 2.2.1}))
\end{align*}
and 
\begin{align*} 
&\dInt((\phi_{\CM}(\ba_{P}^{\bowtie}))^\circ|l_{Q})=\dInt(\bl_P|l_{Q})\\
=&\Int^A(\ba_{P}^{\bowtie}|\gamma_{Q}^{-})+\Int^D(\ba_{P}^{\bowtie}|\gamma_{Q}^{-})+\Int^D(\ba_{P}^{-}|\gamma_{Q}^{-})+ \Int^A(\ba_{P}^{\bowtie}|\gamma_{P}^{\bowtie})+\Int^D(\ba_{P}^{\bowtie}|\gamma_{Q}^{\bowtie})+\Int^D(\ba_{P}^{-}|\gamma_{Q}^{\bowtie})\\
=&\Int^A(\ba_{P}^{\bowtie}|\gamma_{Q}^{-})+\Int^D(\ba_{P}^{\bowtie}|\gamma_{Q}^{-})+ \Int^A(\ba_{P}^{\bowtie}|\gamma_{P}^{\bowtie})+\Int^D(\ba_{P}^{\bowtie}|\gamma_{Q}^{\bowtie}) \ \ \ \ \ \ \ \ \  (\text{by }(\ref{gongshi: D-intersection of a and gammaq in subcase 2.2.1}))\\
=&\Int(\ba_{P}^{\bowtie}|\gamma_{Q}^{-})+\Int(\ba_{P}^{\bowtie}|\gamma_{Q}^{\bowtie})\ \ \ \ \ \ \ \ \ \ \ \ \ \ \ \ \ \ \ \ \ \ \ \ \ \ \ \ \ \ \ \ \ \ \  \ \ \ \ \ \ \ \ \ \ \ \ \ \ \ \ \ (\text{by }(\ref{gongshi: intersection of a and gammaQ in case 2})).
\end{align*}
Now assume that $\sigma_{P}({\CM})< 0$, then $\gamma_Q^-$ and $\gamma_Q^{\bowtie}$ are tagged notched at the end incident to $P$. It follows that
\begin{align}\label{gongshi:D-intersection of a and gamma q in subcase 2.2.2}
    \Int^D(\ba_P^{\bowtie}|\gamma_{Q}^{-})=\Int^D(\ba_P^{\bowtie}|\gamma_{Q}^{\bowtie})=0.
\end{align}
 Hence,  by (\ref{gongshi: intersection of a and gammaQ in case 2}), we have 
\begin{align}\label{gongshi:intersection of ap+ and gammaq in 2.2.2}
  \Int(\ba_{P}^{\bowtie}|\gamma_{Q}^{-})+\Int(\ba_{P}^{\bowtie}|\gamma_{Q}^{\bowtie})=\Int^A(\ba_{P}^{\bowtie}|\gamma_{Q}^{-})+ \Int^A(\ba_{P}^{\bowtie}|\gamma_{Q}^{\bowtie}).  
\end{align}
 Consequently,
\begin{align*} 
 \dInt((\phi_{\CM}(\ba_{P}^{\bowtie}))^\circ|l_{Q})&=\dInt(\ba_{P}^{-}|l_{Q})\\
 &= \Int^A(\ba_{P}^{-}|l_{Q}) &&\\
&=\Int^A(\ba_{P}^{-}|\gamma_{Q}^{-})+\Int^A(\ba_{P}^{-}|\gamma_{Q}^{\bowtie})&& (\text{by (3.1) of Lemma}~\ref{l:relation of a_p and l_Q})\\
&=\Int^A(\ba_{P}^{\bowtie}|\gamma_{Q}^{\bowtie})+\Int^A(\ba_{P}^{\bowtie}|\gamma_{Q}^{\bowtie})&&\\
&=\Int(\ba_{P}^{\bowtie}|\gamma_{Q}^{-})+\Int(\ba_{P}^{\bowtie}|\gamma_{Q}^{\bowtie})&&(\text{by } (\ref{gongshi:intersection of ap+ and gammaq in 2.2.2}))
\end{align*}
 and
 \begin{align*} 
&\dInt((\phi_{\CM}(\ba_{P}^{-}))^\circ|l_{Q})\\=&\dInt(\bl_{P}|l_{Q})\\
=&\Int^A(\ba_{P}^{\bowtie}|\gamma_{Q}^{-})+\Int^D(\ba_{P}^{\bowtie}|\gamma_{Q}^{-})+\Int^D(\ba_{P}^{-}|\gamma_{Q}^{-})+ \Int^A(\ba_{P}^{\bowtie}|\gamma_{Q}^{\bowtie})+\Int^D(\ba_{P}^{\bowtie}|\gamma_{Q}^{\bowtie})+\Int^D(\ba_{P}^{-}|\gamma_{Q}^{\bowtie})\\
=&\Int^A(\ba_{P}^{-}|\gamma_{Q}^{-})+\Int^D(\ba_{P}^{-}|\gamma_{Q}^{-})+ \Int^A(\ba_{P}^{-}|\gamma_{P}^{\bowtie})+\Int^D(\ba_{P}^{-}|\gamma_{Q}^{\bowtie})\ \ \ \ \ \ \ \ \ \ \ \ \ \ \ \ \ \ \  \ \ \ (\text{by } (\ref{gongshi:D-intersection of a and gamma q in subcase 2.2.2}))\\
  =&\Int(\ba_{P}^{-}|\gamma_{Q}^{-})+\Int(\ba_{P}^{-}|\gamma_{Q}^{\bowtie})\ \ \ \ \ \ \ \ \ \ \ \ \ \ \ \ \ \ \ \ \ \ \ \ \ \ \ \ \ \ \ \ \ \ \ \ \ \ \ \ \ \ \ \ \ \ \ \ \ \ \ \ \ \ \ \ \ \ \ \ \ \ \ \  (\text{by } (\ref{gongshi: intersection of a and gammaQ in case 2})).
\end{align*}

\noindent$\mathbf{Case}$ 3: One endpoint of $\ba$ is a puncture ${P}$ and $P=Q$.
Hence $\ba=\ba_P^{-}$ or $\ba_P^{\bowtie}$.  Because $\ba$ is not a loop, we obtain $$\Int^B(\ba|l_{Q})=\Int^B(\ba|\gamma_{Q}^{\bowtie})=\Int^B(\ba|\gamma_{Q}^{-})=0.$$ It is clear that $\Int^D(\ba_P^-|\gamma_Q^{-})=0$,  $\Int^D(\ba_P^{\bowtie}|\gamma_Q^{-})=1$ and  $\Int^D(\ba_P^{-}|\gamma_Q^{\bowtie})=1$. Hence, by (\ref{gongshi:intersection of a and gamma in general case}), we have
$$ \Int(\ba_P^{\bowtie}|\gamma_{Q}^{-})+\Int(\ba_P^{\bowtie}|\gamma_{Q}^{\bowtie})=\Int^A(\ba_P^{\bowtie}|\gamma_{Q}^{-})+\Int^A(\ba_P^{\bowtie}|\gamma_{Q}^{\bowtie})+1,$$   
$$\Int(\ba_P^{-}|\gamma_{Q}^{-})+\Int(\ba_P^{-}|\gamma_{Q}^{\bowtie})=\Int^A(\ba_P^{-}|\gamma_{Q}^{-})+\Int^A(\ba_P^{-}|\gamma_{Q}^{\bowtie})+1.$$
Moreover, 
\begin{align*}
    \dInt(\ba_P^{-}|l_{Q})=&\Int^A(\ba_P^{-}|l_{Q}) && (\text{by (\ref{gongshi: intersection of a and lQ in case 2})} )\\
     =&2\Int^A(\ba_P^{-}|\gamma_{Q}^{-})+1 && (\text{by  (3) of Lemma~\ref{l:relation of a_p and l_Q} })  \\ 
     =&\Int^A(\ba_P^{-}|\gamma_{Q}^{-})+\Int^A(\ba_P^{-}|\gamma_{Q}^{\bowtie})+1
\end{align*}
and
\begin{align*}
\dInt(\bl_{P}|l_{Q})=&\Int^A(\bl_{P}|l_{Q})+\Int^B(\bl_{P}|l_{Q})\\
 =&2\Int^A(\bl_{P}|\gamma_{Q}^{-})+2\Int^B(\bl_{P}|\gamma_{Q}^{-})-1\quad\quad\quad\quad (\text{by (2.2) and (2.3) of Lemma \ref{l:relation of a_p and l_Q}}) \\
   =&2(\Int^A(\ba_P^{\bowtie}|\gamma_{Q}^{-})+\Int^D(\ba_P^{\bowtie}|\gamma_{Q}^{-})+\Int^D(\ba_P^{-}|\gamma_{Q}^{-}))-1 \quad\quad (\text{by  (5) of Lemma~\ref{l:relation of l_p and a_p}})\\
  =&\Int^A(\ba_P^{-}|\gamma_{Q}^{-})+\Int^A(\ba_P^{-}|\gamma_{Q}^{\bowtie})+1.
\end{align*}
Since $\Int^A(\ba_P^{\bowtie}|\gamma_{Q}^{-})=\Int^A(\ba_P^{\bowtie}|\gamma_{Q}^{\bowtie})=\Int^A(\ba_P^{-}|\gamma_{Q}^{-})=\Int^A(\ba_P^{-}|\gamma_{Q}^{\bowtie})$,  we have 
$$\dInt(\ba_P^{-}|l_{Q})=\dInt(\bl_{P}|l_{Q})=\Int(\ba_P^{\bowtie}|\gamma_{Q}^{-})+\Int(\ba_P^{\bowtie}|\gamma_{Q}^{\bowtie})=\Int(\ba_P^{-}|\gamma_{Q}^{-})+\Int(\ba_P^{-}|\gamma_{Q}^{\bowtie}).$$
If $\sigma_{P}({\CM})< 0$, then
\begin{align*}
\dInt((\phi_{\CM}(\ba_P^{\bowtie}))^\circ|l_{Q})=\dInt(\ba_P^{-}|l_{Q})
=\Int(\ba_P^{\bowtie}|\gamma_{Q}^{-})+\Int(\ba_P^{\bowtie}|\gamma_{Q}^{\bowtie}),
\end{align*}
\begin{align*}
\dInt((\phi(\ba_P^{-}))^\circ|l_{Q})=\dInt(\bl_{P}|l_{Q})
=\Int(\ba_P^{-}|\gamma_{Q}^{-})+\Int(\ba_P^{-}|\gamma_{Q}^{\bowtie}).
\end{align*}
Similar discussion yields the result for  $\sigma_{P}({\CM})\ge 0$.
\end{proof}

\subsection{Modification of a multiset}
Let $C$ be a multiset of pairwise compatible tagged arcs. 
We define another multiset $C^{\diamond}$, which consists of plain arcs derived from $C$ through the following process:
\begin{itemize}
    \item [$\mathbf{Step}$ 1:]  replace each conjugate pair $(\gamma_P^{-},\gamma_P^{\bowtie})$ with the corresponding loop $l_P=(\gamma^{\bowtie}_P)^{\circ}$ for every puncture $P$;
    \item [$\mathbf{Step}$ 2:] for each puncture $P$, replace all (notched) tags at $P$ by plain ones.
\end{itemize}
 It is straightforward to see that all arcs in $C^{\diamond}$ are pairwise compatible. It is also clear that $C$ can be recovered from $C^{\diamond}$ and  $\{\sigma_{C}(P)| P\  \text{is a puncture in } \SS\}$ through the following process:
\begin{itemize}
\item [$\mathbf{Step}$ 1:]  tag  notched at the end of each arc which is  incident to some puncture $P$ such that $\sigma_C(P)<0$.
\item  [$\mathbf{Step}$ 2:] for every puncture $P$, replace each once-punctured loop $l$ where $P$ is enclosed  with the corresponding pair  $(\gamma_P^{-},\gamma_P^{\bowtie})$ of conjugate arcs, where   $\gamma_P^{-}$ and $\gamma_P^{\bowtie}$ are tagged in the same way with $l$ at the ends incident to the base point of $l$.
\end{itemize}

\begin{lem} \label{set between tag and untag}
Let $\TT$ be  an admissible tagged triangulation of $(\SS,\MM)$ and $\cm$
  a multiset  of pairwise compatible tagged arcs such that $\cm\cap \TT=\emptyset$. Then under the bijection $\phi_{\cm}$, we have
 $$\Intv_{\TT}(\cm)=\dIntv_{\TT_{\cm}^{\circ}}(\cm^{\diamond}).$$
\end{lem}
\begin{proof}
When $S_{\cm}(P)=\{\mathbf{plain, notched}\}$ for some puncture $P$,
for each pair $(\gamma_P^{-},\gamma_P^{\bowtie})$ in $\cm$, $l_P=(\gamma^{\bowtie}_P)^{\circ}\in \cm^{\diamond}$. 
By Lemma~\ref{l:arc between tag and untag2},
we have $$\dInttv_{\TT}(\gamma_P^{-})+\dInttv_{\TT}(\gamma_P^{\bowtie})=\Intv_{\TT_{\cm}^{\circ}}^{\circ}(l_P).$$ 
Then we can assume that $S_\cm(P)\neq\{\mathbf{plain,notched}\}$   for each puncture $P$.
By Lemma~\ref{l:arc between tag and untag1}, for each $\alpha\in\cm$ we have 
 $$\dInttv_{\TT}(\alpha)=\Intv_{\TT_{\cm}^{\circ}}^{\circ}(\overline{\alpha}).$$
    This completes the proof.
\end{proof}

\subsection{Multisets determined by their intersection vectors}
\begin{thm}\label{t:main theorem-admissible case}
Let $(\SS,\MM)$ be a marked surface with  a strong admissible tagged triangulation $\TT$.  Let $\cm,\cn$ be two multisets consisting of pairwise compatible tagged arcs, such that $\dInttv_{\TT}(\cm)=\dInttv_{\TT}(\cn)$. Then $\cm=\cn$.
\end{thm}
\begin{proof}
Consider $\ba\in\TT$. Suppose that $\ba^{(m)} \subseteq \cm, \ba^{(m+1)}\not\subseteq\cm$ and $\ba^{(n)} \subseteq \cn, \ba^{(n+1)}\not\subseteq\cn$  for some non-negative integers $m,n$. Similar as the proof of Lemma~\ref{l: T has no element in cm and cn}, one can get $m=n$. So $\ba^{({m})}\subset\cm$ if and only if $\ba^{({n})}\subset\cn$.  In the following, we can assume that $\cm\cap \TT=\emptyset=\cn\cap\TT$.

Let $P$ be a puncture in $\SS$, $(\ba_P^-,\ba_P^{\bowtie})$ a pair of conjugate tagged arcs in $\TT$. One can get that $$\sigma_{\cm}(P)=\dIntt^D(\ba_P^{\bowtie}|\cm)-\dIntt^D(\ba_P^{-}|\cm).$$
Then by $$\dIntt(\ba_P^{-}|\cm)=\Int^A(\ba_P^{-}|\cm)+\Int^B(\ba_P^{-}|\cm)+\Int^C(\ba_P^{-}|\cm)+\Int^D(\ba_P^{-}|\cm),$$ $$\dIntt(\ba_P^{\bowtie}|\cm)=\Int^A(\ba_P^{\bowtie}|\cm)+\Int^B(\ba_P^{\bowtie}|\cm)+\Int^C(\ba_P^{\bowtie}|\cm)+\Int^D(\ba_P^{\bowtie}|\cm)$$
and $\Int^*(\ba_P^{-}|\cm)=\Int^*(\ba_P^{\bowtie}|\cm)$, when $*\in\{A,B,C\}$,
we  deduce that
$$\dIntt(\ba_P^{\bowtie}|\cm)-\dIntt(\ba_P^{-}|\cm)=\sigma_{\cm}(P).$$ Similarly, we can obtain  
$$\dIntt(\ba_P^{\bowtie}|\cn)-\dIntt(\ba_P^{-}|\cn)=\sigma_{\cn}(P).$$ Given that $\dInttv_{\TT}(\cm)=\dInttv_{\TT}(\cn)$, we have $\dIntt(\ba_P^{\bowtie}|\cm)=\dIntt(\ba_P^{\bowtie}|\cn)$ and $\dIntt(\ba_P^{-}|\cm)=\dIntt(\ba_P^{-}|\cn)$.
It follows that $$\sigma_{\cm}(P)=\sigma_{\cn}(P)$$ for each puncture $P$. Then we have $\TT_{\cm}=\TT_{\cn}.$ Moreover, we have $\TT_{\cm}^{\circ}=\TT_{\cn}^{\circ}.$

By Lemma~\ref{set between tag and untag}, we have 
  $$\dIntv_{\TT_{\cm}^{\circ}}(\cm^{\diamond})=\dInttv_{\TT}(\cm)=\dInttv_{\TT}(\cn)=\dIntv_{\TT_{\cn}^{\circ}}(\cn^{\diamond})=\dIntv_{\TT_{\cm}^{\circ}}(\cn^{\diamond}).$$
By applying Theorem~\ref{t:case of s has three vertex}, we get $\cm^{\diamond}=\cn^{\diamond}$. Because  $\cm$  is uniquely determined by $\cm^{\diamond}$   and  $\{\sigma_{\cm}(P)|P\text{ is a puncture in } \SS\}$, while  $\cn$ is uniquely determined by $\cn^{\diamond}$   and   $\{\sigma_{\cm}(P)| P\text{ is a puncture in } \SS\}$,  we can conclude that  $\cm=\cn$.
\end{proof}

\section{Denominator conjecture  for some cluster algebras}\label{s:Denominator conjecture  for some cluster algebras}
Let $(\SS,\MM)$ be a marked surface with boundary, $\TT$ be a tagged triangulation of $(\SS,\MM)$. By \cite{FST08}, one can define a cluster algebra $\ma_\TT=\ma(\SS,\MM,\TT)$ with an initial exchange matrix $B_\TT$ defined by $\TT$. According to \cite[Theorem 7.11]{FST08}, 
  there  exists a bijection between the set of all tagged triangulations of  $(\SS,\MM)$ and the set of all clusters of $\ma$, as well as  a bijection between the tagged arcs of $(\SS,\MM)$ and the cluster variables of $\ma$. This bijection maps each tagged arc $\alpha$ to a cluster variable $x[\alpha]$.
Furthermore, as stated in \cite[Thm 8.6]{FST08}, for any two tagged arcs $\alpha$ and $\beta$, the denominator component $d(x[\alpha]|x[\beta])$ in the exchange relation is equal to the intersection number $\Int(\alpha|\beta)$. Specifically, for all tagged arcs $\ba_1,\cdots,\ba_n$ in $\TT$ and any tagged arc $\beta$, we have:
$$x[\beta]=\frac{P(x[\ba_1],\cdots,x[\ba_n])}{\Pi_{i=1}^nx[\ba_i]^{\Int(\ba_i|\beta)}}$$
Here, $P$ is a polynomial that is not divisible by any variable $x[\ba_i]$. This means that the denominator of $x[\beta]$ with respect to the initial seed $(x[\ba_1],x[\ba_2],\cdots,x[\ba_n])$ is $(\Int(\ba_i|\beta))_{\ba_i\in\TT}$, which is equivalent to $\Intv_\TT(\beta)$. 

Recall that a {\it cluster monomial} of $\mathcal{A}_{\TT}$ is a monomial of cluster variables belonging to the same cluster.
Therefore the bijection $\alpha\mapsto x[\alpha]$ yields 
  a bijection $\cm\mapsto x[\cm]$ between  the set of multisets of pairwise compatible tagged arcs of $(\SS,\MM,\TT)$ and the set of cluster monomials in $\ma_T$.
 Clearly, $\dInttv_{\TT}(\cm)=\den(x[\cM])$.

\begin{thm}\label{t:main theorem with strong admissible initial seed}
Let $(\SS,\MM)$ be a marked surface with  is a strong admissible tagged triangulation $\TT$, $\ma_\TT$  the corresponding cluster algebras. Then different cluster monomials in $\ma_\TT$ have different denominators with respect to the initial seed defined by $\TT$.
\end{thm}

\begin{proof}
 Let $\mathbf{m}$ and $\mathbf{n}$ be two cluster monomials in $\ma_\TT$ such that $\den(\mathbf{m})=\den(\mathbf{n})$. 
Let $\cm,\cn$ be two multisets of pairwise compatible tagged arcs such that $x[\cM]=\mathbf{m}$ and $x[\cN]=\mathbf{n}$. Then we have 
$$\dInttv_{\TT}(\cm)=\den(\mathbf{m})=\den(\mathbf{n})=\dInttv_{\TT}(\cn).$$ According to Theorem~\ref{t:main theorem-admissible case}, we conclude that  $\cm=\cn$. Therefore $\mathbf{m}=\mathbf{n}$.
\end{proof}


\begin{thebibliography}{AAAA}




\bibitem[CK06]{CK06}P. Caldero and B. Keller, \emph{From triangulated categories to cluster algebras II}, Ann. Sci. Ecole Norm. Sup. 4eme serie,
\textbf{39} (2006), 983-1009.

\bibitem[CK08]{CK08}P. Caldero and B. Keller,  \emph{From triangulated categories to cluster algebras I}, {Invent. Math.} \textbf{172} (2008), 169--211. 


\bibitem[FeST]{FeST12} A. Felikson, M. Shapiro and P. Tumarkin, \emph{Skew-symmetric cluster algebras of finite mutation type}, J. Eur. Math. Soc. (JEMS) \textbf{14} (2012), 1135-1180.





\bibitem[FoG06]{FoG06} V. Fock and A. Goncharov, \emph{Moduli spaces of local systems and higher Teichmuller theory}, Publ. Math.
 Inst. Hautes Etudes Sci. \textbf{103} (2006), 1-211.
 
 \bibitem[FoG09]{FoG09} V. Fock and A. Goncharov, \emph{Cluster ensembles, quantization and the dilogarithm}, Ann. Sci. Ec. Norm.
 Super. (4) \textbf{42}, (2009), no. 6, 865-930.

\bibitem[FZ02]{FZ02}S. Fomin and A. Zelevinsky, \emph{Cluster algebras I: Foundations}, J. Amer. Math. Soc. \textbf{15}(2)(2002), 497-529.


\bibitem[FZ04]{FZ04}
S. Fomin and A. Zelevinsky, \emph{Cluster algebras: notes for the CDM-03 conference}, in \emph{Current Developments in Mathematics} (International Press, Boston, MA, 2004).

\bibitem[FZ07]{FZ07}
	S. Fomin and A. Zelevinsky, \emph{Cluster algebras IV: Coefficients}, Composito Math. \textbf{143} (2007), 112-164.




\bibitem[FST]{FST08}S. Fomin, M. Shapiro and D. Thurston, \emph{Cluster algebras and triangulated surfaces. I. Cluster complexes}, Acta Math. \textbf{201} (1) (2008) 83-146.


\bibitem[FT]{FT18} S. Fomin and D. Thurston,\emph{Cluster algebras and triangulated surfaces. Part II: Lambda lengths}, Memoirs AMS,
 \textbf{255} (1223), 2018.
 


 
	\bibitem[FuG19]{FG}C. Fu and S. Geng, \emph{On indecomposable $\tau$-rigid modules for cluster-tilted algebras of tame type},  J. Algebra \textbf{531} (2019), 249-282.

 \bibitem[FuG22]{FG2022} C. Fu and S. Geng, \emph{Intersection vectors over tilings with applications to gentle algebras and cluster algebras}, 
arXiv:2212.11497v2.

\bibitem[FuGL]{FGL21a} C. Fu, S. Geng and P. Liu, \emph{Cluster algebras arising from cluster tubes I: integer vectors}, Math. Zeit. \textbf{297} (2021), no. 3-4, 1793-1824.








	

\bibitem[GSV]{GSV} M. Gekhtman, M. Shapiro and A. Vainshtein, \emph{Cluster algebras and Weil-Petersson forms}, Duke Math. J.
 \text{127} (2005), 291-311.

\bibitem[GP]{GP}S. Geng and L. Peng, \emph{The dimension vectors of indecomposable modules of cluster-tilted algebras and
the Fomin-Zelevinsky denominators conjecture},  Acta Math. Sin. (Engl. Ser.) \textbf{28}(3) (2012),  581-586.

\bibitem[GY]{GY20}
Y. Gyoda and T. Yurikusa, \emph{F-matrices of cluster algebras from triangulated surfaces},  Ann. Comb. \textbf{24} (2020), no. 4, 649-695.

 


\bibitem[IL]{IL12} G. Irelli and D. Labardini-Fragoso,\emph{ Quivers with potentials associated to triangulated surfaces, part III: Tagged triangulations and cluster monomials}, Compos. Math. \textbf{148} (2012), 1833-1866.





\bibitem[L]{L09} D. Labardini-Fragoso, \emph{Quivers with potentials associated to triangulated surfaces}, Proc. Lond. Math. Soc. (3) \textbf{98} (3) (2009) 797-839.

\bibitem[M]{M17} M. Mills, \emph{Maximal green sequences for quivers of finite mutation type}, Adv. Math. \textbf{319} (2017), 182-210.

\bibitem[MSW11]{MSW11} G. Musiker, R. Schi er and L. Williams, \emph{Positivity for cluster algebras from surfaces}, Adv. Math. \textbf{227} (2011), 2241-2308.
 
 \bibitem[MSW13]{MSW13} G. Musiker, R. Schi er and L. Williams, \emph{Bases for cluster algebras from surfaces}, Compos. Math. \textbf{149},
 2, (2013) 217-263.
 
 \bibitem[NS]{NS}T. Nakanishi and S. Stella, \emph{Diagrammatic description of $c$-vectors and $d$-vectors of cluster algebras of finite type}, Electron. J. Combin. \textbf{21} (2014), 107 pages.


\bibitem[RS]{RS}
 D. Rupel and S. Stella, \emph{Some consequences of categorification}, SIGMA Symmetry Integrability Geom. Methods Appl. \textbf{16} (2020), Paper No. 007, 8pp. 


 
\bibitem[SZ]{SZ}
P. Sherman and A. Zelevinsky, \emph{Positivity and canonical bases in rank $2$ cluster algebras of finite and affine types}, Moscow Math. J. \textbf{4} (2004), no. 4, 947-974.



\end{thebibliography}
\end{document}